\input amssym.def
\input amssym
\documentstyle[12pt]{article} 
\newcommand{\singlespace}{
    \renewcommand{\baselinestretch}{1}
\large\normalsize}

\singlespace
\begin{document}
\def \Z{\Bbb Z}
\def \R{\Bbb R}
\def \Q{\Bbb Q}
\def \C{\Bbb C}
\def \a{\alpha}
\def \b{\beta}

\newtheorem{th}{Theorem}[section]
\newtheorem{prop}[th]{Proposition}
\newtheorem{coro}[th]{Corollary}
\newtheorem{lemma}[th]{Lemma}
\newtheorem{rem}[th]{Remark}
\makebox[5cm]{} \vspace{1 cm}
\begin{center}
{\LARGE The algebraic structure of relative twisted\\
 vertex operators}

Chongying Dong

Department of Mathematics, University of California, Santa Cruz, CA
95064

James Lepowsky

Department of Mathematics, Rutgers University, New Brunswick,  NJ 08903
\end{center}

\begin{abstract}
Twisted vertex operators based on rational lattices have had many
applications in vertex operator algebra theory and conformal field
theory.  In this paper, ``relativized'' twisted vertex operators are
constructed in a general context based on isometries of rational
lattices, and a generalized twisted Jacobi identity is established for
them.  This result generalizes many previous results.  Relatived
untwisted vertex operators had been studied in a monograph by the
authors.  The present paper includes as a special case the proof of
the main relations among twisted vertex operators based on even
lattices announced some time ago by the second author.
\end{abstract}

\pagestyle{empty}

\section{Introduction}

Vertex operators associated with the roots of a simple Lie algebra
play an important role in the representation theory
of affine Kac-Moody algebras (see for example [BHN], [FK], [Ha], [KKLW], 
[KP], [L1], [LP1]-[LP3], [LW1]-[LW5], [S]).
The algebraic structure of the untwisted vertex 
operators associated with an arbitrary even lattice $L$ is one of the 
important 
motivations for introducing the notion of vertex (operator) algebra [B]
and also provides basic 
examples of  such algebras (see also [FLM3]). Objects that we called
``relative untwisted vertex operators''  
were introduced and studied in [DL1]-[DL2], where the untwisted 
$Z$-algebras ([LP1]-[LP3]), which were used to construct bases of standard
modules for affine Lie algebras (see also [LW2]-[LW5]), were placed
into a systematic axiomatic context, along with 
parafermion algebras [ZF1]. 
Roughly speaking, the relative untwisted vertex operators were defined 
so as to ``relativize'' the usual untwisted vertex operators with respect to
a suitable subspace of the complex span of $L.$ 
The introduction of relative untwisted vertex operators 
led in [DL2] to three levels of generalization of the concept of 
vertex 
operator algebra (and module). One of these notions $-$
 that of ``generalized vertex
operator algebra'' $-$ enabled us to clarify  the essential equivalence
between the $Z$-algebras of [LP1]-[LP2] and the 
parafermion algebras of [ZF1], among
other things providing a precise mathematical foundation for ``parafermion conformal field theory.''

\pagestyle{plain}
General {\it twisted}
 vertex operators associated with an arbitrary even lattice 
$L$ equipped
with
a finite order isometry were introduced and studied in [L1] and [FLM2];
see also [KP] for a different approach in  the special case of the elements 
of root lattices corresponding to the roots of a simply-laced simple Lie algebra (this approach does not generalize to arbitrary root-lattice elements,
which are treated in [L1]).
The twisted vertex operators, parametrized by
an ``untwisted space'' $V_L$ and acting on a ``twisted space'' $V_L^T,$ satisfy
a ``twisted Jacobi identity'' (see [FLM2]-[FLM3] and [L1]-[L2]). 
The important and suggestive case in which the isometry is $-1$ 
plays a special role in constructing the moonshine module vertex operator 
algebra $V^{\natural}$ (see [FLM1]-[FLM3]).
Twisted vertex operators are closely related to orbifold theory (see 
[DVVV], [DHVW] and [DM2]).
 
In this  paper, we present a general theory of twisted
vertex operators relativized with respect to a suitable subspace of 
the complex 
span of the rational lattice $L,$ and based on lattice 
automorphisms, generalizing
and incorporating the construction of and the results on the ordinary
twisted vertex operators announced in [L1] and [FLM2] (see also
[L2]). The main result of the present paper is 
a generalized 
twisted Jacobi identity for relative twisted vertex operators. 
This single result generalizes many known ones.
At the same time we supply detailed proofs of the main 
results announced in [L1]
and [FLM2]. These results, for the case of ordinary (not relative) 
twisted vertex operators, 
were the main motivation for introducing the notion of twisted modules
for a vertex algebra, in [FFR] and [D]. In particular,
the twisted spaces $V_L^{T}$ constructed in [L1] and [FLM1]-[FLM3] are  
twisted modules for the vertex algebras $V_L$ associated to 
even lattices $L.$ 
Recently, our result on (ordinary) 
twisted vertex operators (for non-even lattices) has been  
used in [DM1] to give another construction of the moonshine module
associated with a certain isometry  (of order 3) of the Leech lattice. 
Also, twisted vertex operators based on an integral lattice
have been studied in [X].

The special case in which $L$ is the weight lattice of $sl(2,{\Bbb C})$ and the
lattice isometry is $-1$ has been studied previously by Husu [Hu],
where the generalized twisted Jacobi identity was established
for relative twisted vertex operators. As an application of this
identity, many results concerning the twisted $Z$-algebras discovered
and developed in 
[LW2]-[LW4] in the course of constructing  bases of the standard
modules for the affine Lie algebra $A_1^{(1)}$ were reinterpreted
in a natural way from the vertex-operator-theoretic viewpoint. As
pointed
out in [Hu], the constructions in [LW2]-[LW4] and [Hu] 
correspond to the (twisted)
parafermion algebras of [ZF2].

This paper is organized as follows: In Section 2 we present a basic
setting, involving a rational lattice $L$ and the vector space which it
generates,
and an isometry of finite order of the lattice. We also consider
a subspace of the vector space with respect to which the vertex
operators based on  the untwisted space $V_L$
will
be ``relativized.'' In Section 3, we recall the notion of relative 
untwisted vertex operator (see [DL2]). These operators
are parametrized by the space $V_L$
and also act on this space. Relative
twisted vertex operators, parametrized by the untwisted space and
acting on the twisted space $V_L^T,$ 
are defined in Section 4. A fundamental quadratic operator $\Delta_z$
enters into the construction of these operators, as in [FLM2] and [FLM3].    

Section 5, the central part of this work, is devoted to the
formulation and proof of our generalized twisted Jacobi identity for
relative twisted vertex operators. In Section 6, we focus our
attention on certain relative (untwisted and twisted) vertex operators
providing representations of the Virasoro algebra, and we also compute
the weights of the twisted spaces.  The quadratic operator $\Delta_z$
plays a basic role here. We observe that the resulting shift in
conformal weight associated with the twisted space is related to
certain values of the second Bernoulli polynomial.  Finally, in
Section 7, we put the twisted spaces $V_L^T$ constructed in [L1] and
[FLM2] (see also [L2]) into the axiomatic context of twisted modules
for the vertex algebra $V_L.$

{\it Acknowledgments.} During the course of this work, much of 
which was done many years ago, C.D. was partially
supported by a postdoctoral fellowship at Rutgers University,
a Regents' Junior Faculty Fellowship from the 
University of California, faculty research funds granted by the University 
of California at Santa Cruz, NSA grant  
MDA904-92-H-3099 and NSF grant DMS-9303374; and J.L. by a 
Guggenheim Foundation Fellowship, the Rutgers University Faculty
Academic Study Program and NSF grants DMS-8603151 and DMS-9111945;
see also the acknowledgments in the announcement [L1], the full 
proofs of whose main results on relations among twisted 
vertex operators are finally supplied in the present paper.

\section{The setting}

In this section we introduce notation and  assumptions. 
In particular, we introduce a 
lattice $L$ with an
isometry $\nu,$ two central extensions $\hat L$ and  
$\hat L_{\nu}$ by the same finite cyclic
group, and a subspace ${\bf h}_*$ of the complex span ${\bf h}$ of
$L.$ The lattice $L$ together with the central extension $\hat L$ will be used
to build a vector space $V_L$ which will parametrize  
relative untwisted vertex operators and on which these operators will
also act. The space $V_L$ will also parametrize relative twisted
vertex operators acting on a vector space $V_L^T$ that will be built from
$L$ and the central extension $\hat L_{\nu}.$ 

We work in the setting of [L1], [FLM2] and [DL2].
Let $p$ and $q$ be two fixed positive integers
such that $p$ divides $q.$  
The following data and conditions are 
assumed:

2.1. Let $L$\  be a (rational) lattice equipped  with a symmetric
nondegenerate ${\Bbb Q}$-valued ${\Bbb Z}$-bilinear form\ $ \langle
\cdot,\cdot\rangle ,$ not necessarily positive definite. 

2.2. Let $\nu$\  
be an isometry of\  $L$\  such that\  $\nu^p=1$ (but $p$
need not be the exact order of $\nu$ and we may take $\nu=1$).

2.3. Let $c_0$ and $c_0^{\nu}$ be two $\nu$-invariant  alternating 
${\Bbb Z}$-bilinear maps
\begin{equation}
L\times L\longrightarrow {\Bbb Z}/q{\Bbb Z}.
\end{equation}

\begin{rem}\label{rem1.1}\hspace{-0.2 cm}{\bf :\ }{\rm  For a positive 
integer $n$ let $\langle\kappa_n\rangle$ be the cyclic group 
of order $n$ with generator $\kappa_n.$ We assume that
the generators are related by the condition
$\kappa_n=\kappa_m^{m/n}$ if $n$ divides $m.$ (This holds if
for instance 
$\kappa_n=e^{2\pi i/n}$ for all $n.$)
Viewed as alternating ${\Bbb Z}$-bilinear maps from $L\times L$
to ${\Bbb Z}/q{\Bbb Z},$ $c_0$ and $c_0^{\nu}$ determine two  central
extensions of $L$,
\begin{equation}
\begin{array}{c}
1\ \rightarrow \ \langle \kappa_q\rangle \ \rightarrow \ \hat{L} \
\bar{\rightarrow} \ L \ \rightarrow \ 1\\
1\ \rightarrow \ \langle \kappa_q\rangle \ \rightarrow \ \hat{L}_{\nu}
\ \bar{\rightarrow} \ L \ \rightarrow \ 1,
\end{array}
\end{equation}
uniquely up to equivalence, by the commutator conditions
\begin{equation}
\begin{array}{c}
aba^{-1}b^{-1}=\kappa_q^{c_0(\bar{a},\bar{b})},\ \ aba^{-1}b^{-1}=
\kappa_q^{c_0^{\nu}(\bar{a},\bar{b})}
\in \langle \kappa_q\rangle 
\end{array}
\end{equation}
for $a, b\in \hat{L}$ or $a,b\in\hat{L}_{\nu}$ (cf. for example [FLM3],
Section 5.2).
Then there is a  
set-theoretic identification (which is usually not an isomorphism of
of groups)  between the groups $\hat{L}$ and $\hat{L}_{\nu}$ such that
the respective group multiplications $\times$ and $\times_{\nu}$
are related by
\begin{equation}\label{2.6}
a\times b=\kappa_q^{\epsilon_0(\bar a,\bar b)}a\times_{\nu}b,
\end{equation}
where $\epsilon_0$ is a bilinear map from $L\times L$ to ${\Bbb
Z}/q{\Bbb Z}$ satisfying 
\begin{equation}
\epsilon_0(\alpha,\beta)-\epsilon_0(\beta,\alpha)=c_0(\alpha,\beta)
-c_0^{\nu}(\alpha,\beta)
\end{equation}
for $\alpha,\beta\in L.$ (To see this, we may choose a section
of $\hat L$ and a section of $\hat L_{\nu}$ such that
the corresponding 2-cocycles are bilinear and we then take 
$\epsilon_0$ to be the difference of the two cocycles. The identification
is evident.)
Moreover, $\nu$ lifts to an automorphism
$\hat{\nu}$ of $\hat{L}$ 
fixing $\kappa_q$, that is, 
\begin{equation}
\overline{\hat{\nu}a}=\nu\bar{a} \ \ \ \ \ {\rm for}\ \ a\in \hat{L},\
\ \ \ \ \hat{\nu}\kappa_q=\kappa_q,
\end{equation}
and such a lifting is unique up to multiplication by a lifting of the
identity automorphism of $L$, which is of the form 
\begin{equation}
\begin{array}{r}
\lambda^*: \hat{L}\to \hat{L}\ \ \ \  \ \\
 a\mapsto a\kappa_q^{\lambda(\bar{a})} 
\end{array}
\end{equation}
for some $\lambda\in$Hom$(L,{\Bbb Z}/q{\Bbb Z})$ (cf. [FLM3], Section
5.4). If $\nu=1,$ then we may take $\hat{\nu}=1.$ Note that $\hat\nu$ also
acts set-theoretically on $\hat L_{\nu}$ by using the identification. 
Moreover, if 
$\epsilon_0(\cdot,\cdot)$ is $\nu$-invariant, then $\hat\nu$ 
is also an automorphism of $\hat{L}_{\nu},$ from (\ref{2.6}).}
\end{rem}
\begin{rem}\label{rem1.2}\hspace{-0.2 cm}{\bf :\ }{\rm  Here we describe
the special case which is the subject of [L1]; see also Sections 3.1 and 
3.3 of [FLM2]. Let $L$ be
an even lattice, i.e., $\langle \alpha,\alpha\rangle \in 2{\Bbb Z}$ for all
$\alpha\in L.$ Take $q=p$ if $p$ is even and take $q=2p$ if $p$ is odd, so
that $q$ is always even.
We assume that 
\begin{equation}\label{L2.1}
\langle\nu^{p/2}\alpha,\alpha\rangle\in 2{\Bbb Z}\ \ \ {\rm for}\ \ \
\alpha\in L
\end{equation}
if $p$ is even; this can always be arranged by doubling $p$ if necessary. Then
\begin{equation}\label{L2.2}
c_0(\alpha,\beta)=q\langle \alpha,\beta\rangle/2 +q{\Bbb Z},
\end{equation}
\begin{equation}\label{add1}
c^{\nu}_0(\alpha,\beta)={\displaystyle \sum_{r=0}^{p-1}(q/2+qr/p)\langle\
\nu^r\alpha,\beta\rangle +q{\Bbb Z}}
\end{equation} 
for $\alpha,\beta\in L$  define two $\nu$-invariant 
alternating ${\Bbb Z}$-bilinear 
maps from $L\times L$ to ${\Bbb Z}/q{\Bbb Z}.$
 We have the commutator relations  
\begin{equation}
aba^{-1}b^{-1}=\kappa_2^{\langle\bar a,\bar b\rangle} 
\end{equation}
in $\hat L$ and
\begin{equation}
aba^{-1}b^{-1}=\kappa_2^{\sum_{r=0}^{p-1}\langle\nu^r\bar a,\bar b\rangle} 
\kappa_p^{\sum_{r=0}^{p-1}r\langle\nu^r\bar a,\bar b\rangle}
\end{equation}
in $\hat L_{\nu}.$
We may take 
\begin{equation}\label{LL.2}
\epsilon_0(\alpha,\beta)=\sum_{0<r<p/2}(q/2+qr/p)\langle\nu^{-r}\alpha,
\beta\rangle+q{\Bbb Z}
\end{equation}
in this case. Note that $\epsilon_0(\cdot,\cdot)$ 
is $\nu$-invariant  and thus $\hat\nu$ is also an automorphism 
of $\hat L_{\nu}.$ It is shown in [L1] that there is a lifting 
$\hat\nu$ of $\nu$ with the special property that for $a\in\hat L,$ $\hat{\nu}a=a$ if 
$\nu\bar{a}=\bar{a}.$ In the special cases $p=1$ and $p=2,$ 
we have $c_0^{\nu}=c_0, \epsilon_0=0,$ and $\hat L_{\nu}=\hat L$ as groups.}
\end{rem}

2.4. Let ${\bf h}_*$ be a subspace (possibly $0$ or ${\bf h}$) of 
\begin{equation}
{\bf h}={\Bbb C}\otimes_{\Bbb Z}L
\end{equation}
on which the natural (nonsingular) extension of the form $\langle \cdot,\cdot\rangle $ on $L$,
still denoted $\langle \cdot,\cdot\rangle $, remains nonsingular. That is, 
\begin{equation}\label{2.1}
{\bf h}={\bf h}_*\oplus {\bf h}^{\bot}_*,
\end{equation}
${}^{\bot}$ denoting orthogonal complement. We write
\begin{equation}
\begin{array}{l}
{\bf h}\to {\bf h}^{\bot}_*\ \ \ \ \ \ \ \ {\bf h}\to {\bf h}_*\\
h\mapsto h', \ \ \ \ \ \ \ \   h\mapsto h''
\end{array}
\end{equation}
for the projection maps to ${\bf h}^{\bot}_*$ and ${\bf h}_*.$ We also
assume that ${\bf h}_*$ is stable under the natural action of $\nu$ on
${\bf h}:$
\begin{equation}
\nu {\bf h}_*={\bf h}_*.
\end{equation}
Then the two projection maps commute with the action of $\nu.$

\section{Relative untwisted vertex operators}
\setcounter{equation}{0}
\setcounter{th}{0}

We shall define the untwisted space $V_L$ (cf. [FLM3]) and 
present the notion 
 of relative untwisted 
vertex operators $Y_*(v,z)$ ($v\in V_L$) acting on $V_L,$ following [DL1] and  [DL2]. This section, almost a copy of Chapter 3 of [DL2], is supplied for the
reader's convenience. The degenerate case ${\bf h}_*=0$ (the 
``unrelativized'' case) is important in its
own right.

The affine Lie algebra associated with the abelian Lie algebra
${\bf h}$ is given by:
\begin{equation}\label{3.1}
\hat{\bf h}={\bf h}\otimes{\Bbb C}[t,t^{-1}]\oplus
{\Bbb C}c
\end{equation}
with structure defined by
\begin{equation}[x\otimes t^m, y\otimes t^n]=\langle x,y\rangle m\delta_{m+n,0}c\ \
\mbox{for} \ \  x,y\in{\bf h},\ \ m,n\in{{\Bbb Z}}, 
\end{equation}
\begin{equation}\label{3.3}
[c,\hat{\bf h}]=0.
\end{equation}
(This notation and related notation below may be applied to any 
finite-dimensional abelian Lie algebra with a nonsingular symmetric
form.) By (\ref{2.1}) we have
\begin{equation}\label{3.4}
\hat{\bf h}={\bf h}_*\otimes{\Bbb C}[t,t^{-1}]\oplus{\bf h}_*^{\perp}
\otimes{\Bbb C}[t,t^{-1}]\oplus{\Bbb C}c.
\end{equation}
Then $\hat{\bf h}$ has a ${\Bbb Z}$-gradation, the {\it weight gradation}
associated with ${\bf h}_*,$ given by:
\begin{equation}\label{3.5}
{\rm wt}\,(x\otimes t^m)=0,\ \ \ {\rm wt}\,(y\otimes t^n)=-n,\ \ \
{\rm wt}\,
c=0
\end{equation}
for $x\in {\bf h}_*$, $y\in{\bf h}_*^{\perp}$ and $m,n\in {\Bbb Z}.$

Set
\begin{equation}\label{3.6}
\hat{\bf h}^+={\bf h}\otimes t{\Bbb C}[t],\ \  \hat{\bf h}^-={\bf
h}\otimes t^{-1}{\Bbb C}[t^{-1}]. 
\end{equation}
The subalgebra                                
\begin{equation}\label{3.7}
\hat{\bf h}_{\Bbb Z}= 
\hat{\bf h}^+\oplus\hat{\bf h}^-\oplus {\Bbb C}c 
\end{equation}
of $\hat{\bf h}$ is a Heisenberg algebra. Consider the
induced $\hat{\bf h}$-module, irreducible even under $\hat{\bf h}_{\Bbb Z},$ 
\begin{equation}\label{3.8}
M(1)=U(\hat{\bf h})\otimes_{U({\bf h}\otimes{{\Bbb C}}[t]
\oplus{{\Bbb C}}c)}{{\Bbb C}}\simeq S(\hat{\bf h}^-)\ \ \ (\mbox{linearly}),
\end{equation}
${\bf h}\otimes{{\Bbb C}}[t]$ acting trivially on ${{\Bbb C}}$ and $c$
acting as 1; $U(\cdot)$ denotes universal enveloping algebra and 
$S(\cdot)$ denotes symmetric algebra.  The $\hat{\bf
h}$-module $M(1)$ is ${\Bbb Z}$-graded so that wt\,1\ =\ 0 (we write 1 for
$1\otimes 1$):
\begin{equation}
M(1)=\coprod_{n\in {\Bbb Z},\,n\ge 0}M(1)_n,
\end{equation}
where $M(1)_n$ denotes the homogeneous subspace of weight $n.$
The automorphism
$\nu$ of $L$ acts in a natural way on ${\bf h},$ on $\hat{\bf h}$
(fixing $c$) and on $M(1)$, preserving the gradations, and for $u\in \hat{\bf h}$ and $m\in
M(1),$
\begin{equation}\label{e3.10}
\nu (u\cdot m)=\nu (u)\cdot \nu(m).
\end{equation}
 
Form the induced $\hat{L}$-module and ${\Bbb C}$-algebra 
\begin{equation}\label{3.11a}
\begin{array}{l}
{{\Bbb C}}\{L\}={\Bbb C}[\hat{L}]/(\kappa_q-\omega_q){\Bbb C}[\hat{L}]\\
\ \ \ \ \ \ \ \ \!\ ={\Bbb C}[\hat{L}]\otimes_{{\Bbb C}[\langle \kappa_q\rangle ]}{\Bbb C}\simeq{{\Bbb
C}}[L]\ \ (\mbox{linearly}),
\end{array} 
\end{equation}
where ${{\Bbb C}}[\cdot]$ denotes group algebra, $\omega_q$ is a
primitive $q^{th}$ root of unity in ${\Bbb C}^{\times}$ and 
$\kappa_q$ acts as $\omega_q$ on ${\Bbb C}.$  For $a\in\hat{L}$, 
write $\iota(a)$ for the image of $a$ in ${{\Bbb C}}\{L\}$. Then the action of 
$\hat{L}$ on ${{\Bbb C}}\{L\}$ and the product in ${\Bbb C}\{L\}$ are given by: 
\begin{equation}
a\cdot\iota(b)=\iota(a)\iota(b)=\iota(ab),
\end{equation}
\begin{equation}
\kappa_q\cdot\iota(b)=\omega_q\iota(b)
\end{equation}
for $a,b\in\hat{L}.$ We give ${\Bbb C}\{L\}$ the ${\Bbb C}$-gradation
determined by:
\begin{equation}
\mbox{wt}\,\iota(a)=\frac{1}{2}\langle \bar{a}',\bar{a}'\rangle  \ \ \ \ \mbox{for}\ \
a\in \hat{L}.
\end{equation}
The automorphism $\hat{\nu}$ of $\hat{L}$  acts canonically on ${\Bbb C}\{L\},$ preserving the gradation, 
in such a way that 
\begin{equation}\label{e3.15}
\hat{\nu}\iota(a)=\iota(\hat{\nu}a) \ \ \ \ \mbox{for}\ \ \
a\in\hat{L},
\end{equation}
and we have 
\begin{equation}
\hat{\nu}(\iota(a)\iota(b))=\hat{\nu}(a\cdot \iota(b))=\hat{\nu}(a)\cdot\hat{\nu}\iota(b)=\hat{\nu}\iota(a)\hat{\nu}\iota(b).
\end{equation}

Also define a grading-preserving action of ${\bf h}$ on ${\Bbb
C}\{L\}$ by:
\begin{equation}
h\cdot\iota(a)=
\langle h',\bar{a}\rangle \iota(a)
\end{equation}
for $h\in{\bf h}$ (so that ${\bf h}_*$ acts trivially). Then ${\bf h}$
acts as algebra derivations and
\begin{equation}
\hat{\nu}(h\cdot \iota(a))=\nu(h)\cdot \hat{\nu}\iota(a).
\end{equation}
We shall use a formal variable $z$ (and later,  commuting formal
variables $z, z_0, z_1, z_2,$ etc.). Define 
\begin{equation}\label{3.19}
z^h\cdot\iota(a)=z^{\langle h',\bar{a}\rangle }\iota(a) 
\end{equation}
for $h\in{\bf h}.$ Then
\begin{equation}
\hat{\nu}(z^h\cdot\iota(a))=z^{\nu(h)}\cdot\hat{\nu}\iota(a).
\end{equation}
We shall mostly be interested in the actions of  $h$ and $z^h$ on ${\Bbb
C}\{L\}$ only for $h\in{\bf h}^{\perp}_*.$

Set
\begin{equation}\label{3.21}
V_L=M(1)\otimes_{\Bbb C}{\Bbb C}\{L\}\simeq S(\hat{\bf h}^-)\otimes {\Bbb
C}[L] \ \ \ \ (\mbox{linearly})
\end{equation} 
and give $V_L$ the tensor product ${\Bbb C}$-gradation:
\begin{equation}
V_L=\coprod_{n\in {\Bbb C}}(V_L)_n.
\end{equation}
We have wt$\,\iota(1)=0,$ where we identify ${\Bbb C}\{L\}$ with $1\otimes {\Bbb C}\{L\}.$
Then $\hat{L},$
$\hat{\bf h}_{\Bbb Z},$ $h,$ $z^h$ $(h\in{\bf h})$ act naturally on
$V_L$ by acting on either $M(1)$ or
${{\Bbb C}}\{L\}$ as indicated above. In particular, $c$ acts as 1 and
${\bf h}_*$ acts trivially. The
automorphism $\hat{\nu}$ acts in a natural grading-preserving way on $V_L,$ via
$\nu\otimes\hat{\nu},$ and this action is compatible with the other
actions:
\begin{equation}\label{3.23}
\hat{\nu}(a\cdot v)=\hat{\nu}(a)\cdot\hat{\nu}(v)
\end{equation}
\begin{equation}\label{3.24}
\hat{\nu}(u\cdot v)=\nu(u)\cdot \hat{\nu}(v)
\end{equation}
\begin{equation}\label{3.25}
\hat{\nu}(z^h\cdot v)=z^{\nu(h)}\cdot \hat{\nu}(v)
\end{equation}
for $a\in \hat{L}, u\in \hat{\bf h}, h\in {\bf h}, v\in V_L.$

It will be convenient later to have the notation
\begin{equation}
c(\alpha,\beta)=\omega_q^{c_0(\alpha,\beta)}\in {\Bbb C}^{\times}
\end{equation}
for all $\alpha, \beta\in L.$ Then
\begin{equation}
aba^{-1}b^{-1}=c(\bar{a},\bar{b})
\end{equation}
for $a,b\in \hat{L},$ as operators on ${\Bbb C}\{L\}$ and on $V_L$, and
\begin{equation}
c(\nu\alpha,\nu\beta)=c(\alpha,\beta).
\end{equation}

For $\alpha\in{\bf h},$
$n\in{{\Bbb Z}}$,
we write $\alpha(n)$ for the operator on $V_L$ determined by $\alpha\otimes
t^n$. Set 
\begin{equation}\label{3.29}
\Omega_\ast=\{v\in V_L \mid h(n)v=0 \ \ \mbox{for}\ \ h\in {\bf
h}_\ast,\ \ n>0 \},
\end{equation}
\begin{equation}\label{3.30}
V_*=\mbox{span}\{h(n)V_L\mid h\in {\bf h}_*,\ \ n<0\}.
\end{equation}
Then $\Omega_*$ is the vacuum space for the Heisenberg algebra 
$(\hat{\bf h}_*)_{\Bbb Z}$ (defined as in (\ref{3.7}))
and we have
\begin{equation}\label{3.31}
V_L=\Omega_*\oplus V_*. 
\end{equation}
In fact, we see from (\ref{3.21}) that $V_L$ has the
decomposition
\begin{equation}
V_L=S(\hat{\bf h}^-_*)\otimes S((\hat{\bf h}_*^{\perp})^-)\otimes
{\Bbb C}\{L\}
\end{equation}
and so
\begin{equation}\label{3.33}
\Omega_*=S((\hat{\bf h}_*^{\perp})^-)\otimes{\Bbb C}\{L\},
\end{equation}
\begin{equation}\label{3.34}
V_*=\hat{\bf h}^-_*S(\hat{\bf h}^-_*)\otimes\Omega_*.
\end{equation}
Here $\hat{\bf h}^-_*$ and $(\hat{\bf h}_*^{\perp})^-$ are defined as
in (\ref{3.6}).
In terms of the general structure of modules for Heisenberg algebras,
we know from [FLM3], Section 1.7, that the (well-defined) canonical
linear map
\begin{equation}
\begin{array}{lcr}
U((\hat{\bf h}_*)_{\Bbb Z})\otimes_{(\hat{\bf h}_*^+\oplus{\Bbb
C}c)}\Omega_*&\to& V_L \\ 
\hspace*{2.8 cm}u\otimes v &\mapsto& u\cdot v 
\end{array}
\end{equation}
$(u\in U((\hat{\bf h}_*)_{\Bbb Z}), v\in \Omega_*)$ is an $(\hat{\bf
h}_*)_{\Bbb Z}$-module isomorphism, and in particular, that the linear
map
\begin{equation}
\begin{array}{lcr}
M_*(1)\otimes_{\Bbb C} \Omega_*=U(\hat{\bf h}_*^-)\otimes_{\Bbb C}\Omega_*&\to& V_L \\ 
\hspace*{4 cm}u\otimes v &\mapsto& u\cdot v 
\end{array}
\end{equation}
$(u\in U(\hat{\bf h}_*^-), v\in \Omega_*)$ is an $(\hat{\bf
h}_*)_{\Bbb Z}$-module isomorphism, $\Omega_*$ now being regarded as a
trivial $(\hat{\bf h}_*)_{\Bbb Z}$-module, where
$M_*(1)$ is the $(\hat{\bf h}_*)_{\Bbb Z}$-module defined by analogy
with (\ref{3.8}). The spaces
$\Omega_*$ and $V_*$ are ${\Bbb C}$-graded and are stable under the
actions of $\hat{\bf h}^{\perp}_*$ (defined as in
(\ref{3.1})), of $\hat{L}$ and of $\hat{\nu}.$

For $\alpha\in{\bf h},$ set
\begin{equation}
\alpha(z)=\sum_{n\in{{\Bbb Z}}}\alpha(n)z^{-n-1}.\end{equation} 
We use a normal ordering procedure, indicated by open colons, which signify
that the enclosed expression is to be reordered if necessary so that
all the operators
$\alpha(n)$ $(\alpha\in{\bf h},\  n<0),$ $a\in \hat{L}$ are to be placed to
the left of all the operators $\alpha(n),$ $z^{\alpha}$ $
(\alpha\in {\bf h},\ n\ge 0)$ before the expression is evaluated. For $a \in \hat{L}$, set 
\begin{equation}\label{3.38}
Y_*(a,z)=\mbox{$\circ\atop\circ$}e^{\int(\bar{a}'(z)-
\bar{a}'(0)z^{-1})}az^{\bar{a}'}\mbox{$\circ\atop\circ$},
\end{equation}
using an obvious formal integration notation. (Note that the symbol 
$z^{\bar{a}'}$ could be replaced by $z^{\bar{a}}$ in this formula, in
view of (\ref{3.19}).) Let
$$a \in \hat{L},\ \alpha_1, ..., \alpha_k \in
{\bf h},\  n_1, ..., n_k \in{\Bbb Z}\ \ (n_i>0)$$
and set
\begin{equation}\label{in}
\begin{array}{l}
v=\alpha_1(-n_1)\cdot
\cdot\cdot\alpha_k(-n_k)\otimes\iota(a)\\
\ \,\ =\alpha_1(-n_1)\cdot
\cdot\cdot\alpha_k(-n_k)\cdot\iota(a)\in V_L. 
\end{array}
\end{equation}
We define
\begin{equation}\label{3.40}
Y_*(v,z)=\mbox{$\circ\atop\circ$}\left(\frac{1}{(n_1-1)!}\left(\frac{
d}{dz}\right)^{n_1-1}\alpha'_1(z)\right)\cdot\cdot\cdot\left(\frac{1}{(n_k-1)!}\left(\frac{d}{dz}\right)^{n_k-1}\alpha'_k(z)\right)Y_*(a,z)\mbox{$\circ\atop\circ$}.
\end{equation}
This gives us a well-defined linear map
\begin{equation}\label{3.41}
\begin{array}{lcr}
V_L&\rightarrow&(\mbox{End}\,V_L)\{z\} \hspace*{4.5 cm}\\
\ v&\mapsto& Y_*(v,z)=\displaystyle{\sum_{n\in{\Bbb C}}}v_nz^{-n-1}, \ \ v_n\in\mbox{End}\,V_L,
\end{array}
\end{equation}
where for any vector space $W$, we define $W\{z\}$ to be the vector
space of $W$-valued formal series in
$z$, with arbitrary complex powers of $z$ allowed:
\begin{equation}
W\{z\}=\{\sum_{n\in {\Bbb C}}w_nz^n\ |\ w_n\in W\}. 
\end{equation}
We call  $Y_*(v,z)$ the $untwisted$ $vertex$ $operator$ $associated$ $with$
$v,$ $defined$ $relative$ $to$ ${\bf h}_{\ast}$. We shall be
especially interested in these operators for $v\in \Omega_*.$ Note
that the component operators $v_n$ of $Y_*(v,z)$ are defined in (\ref{3.41}).

\noindent {\bf Remark 3.2:} The case ${\bf h}_*=0$ recovers the case
of ordinary untwisted vertex operators as defined in [B] and
[FLM3]. In this case, $\Omega_*=V_L,$ $V_*=0,$ and
$Y_*$ is the operator $Y,$ in the notation of [FLM3].
 
It is 
easy to check
{}from the definition that the relative untwisted  vertex operators
$Y_*(v,z)$  for $v\in V_L$  have the following properties:
\begin{equation}
Y_*(1,z)=1
\end{equation}
\begin{equation}
Y_*(v,z)=0 \ \ \ \mbox{if}\ \  v\in V_*
\end{equation}
\begin{equation}\label{per1}
Y_*(v,z)\iota(1)\in V_L[[z]] \ \ \mbox{and}\ \ \ \lim_{z\to
0}Y_*(v,z)\iota(1)=v\ \ \ \mbox{if}  \  v\in \Omega_*
\end{equation}
\begin{equation}\label{3.46}
[(\hat{\bf h}_*)_{\Bbb Z},Y_*(v,z)]=0
\end{equation}
\begin{equation}
Y_*(a\cdot v,z)=aY_*(v,z)
\end{equation}
for $a\in\hat{L}$ such that $\bar{a}\in L\cap {\bf h}_*.$ (For a
vector space $W,$ $W[[z]]$ signifies the space of formal power series in $z$
with coefficients in $W.$) Moreover, using (\ref{3.23})-(\ref{3.25}) we see that 
\begin{equation}\label{e3.48}
\hat\nu Y_*(v,z)\hat\nu^{-1}=Y_*(\hat\nu v,z).
\end{equation}

By (\ref{3.46}), the component operators $v_n$ of $Y_*(v,z)$  preserve 
$\Omega_*$  and $V_*$, so that we get a well-defined linear map
\begin{equation}
\begin{array}{lcr}
V_L&\rightarrow&(\mbox{End}\,\Omega_*)\{z\} \\
\ v&\mapsto& Y_*(v,z).\hspace*{0.5 cm}
\end{array}
\end{equation}
We continue to denote the restriction of this map to $\Omega_*$ by $Y_*:$
\begin{equation}
\begin{array}{lcr}
\Omega_*&\rightarrow&(\mbox{End}\,\Omega_*)\{z\} \\
\ v&\mapsto& Y_*(v,z).\hspace*{0.5 cm} 
\end{array}
\end{equation}
A basic property of the operators (\ref{3.38}) and
ultimately, of the operators (\ref{3.40}), is: For $a,b\in \hat{L},$
\begin{equation}\label{3.51}
Y_*(a,z_1)Y_*(b,z_2)=\mbox{$\circ\atop\circ$}Y_*(a,z_1)Y_*(b,z_2)
\mbox{$\circ\atop\circ$}(z_1-z_2)^{\langle \bar{a}',\bar{b}'\rangle },
\end{equation}
where the binomial expression is to be expanded in nonnegative integral
powers of the second variable, $z_2.$

\noindent{\bf Remark 3.3: } In our frequent use of
formal series, it will typically be understood that a binomial
expression such as that in (\ref{3.51}) is to be expanded as in (\ref{3.51}) $-$
as a formal power series in the second variable.

\section{Relative twisted vertex operators}
\setcounter{equation}{0}
In this section, we systematically construct a setting which
generalizes that of Section 3 to incorporate the lattice-isometry
$\nu.$  Many definitions and notations in Section 3 have natural 
generalizations in this setting. We shall define a twisted space $V_L^T$ and
present a notion of relative twisted vertex operator 
$Y_*^{\nu}(v,z)$ $(v\in V_L)$ acting on $V_L^T.$  We shall point out 
that the twisted vertex operators 
introduced in [L1] and [FLM2] are special cases of our setting.

For $n\in {\Bbb Z},$ set
\begin{equation}\label{9.1}
{\bf h}_{(n)}=\{\alpha \in {\bf h}|\nu
\alpha=\omega_p^n\alpha\}\subset{\bf h}
\end{equation}
where $\omega_p=\omega_q^{q/p},$ which 
is a primitive $p^{th}$ root of unity in ${\Bbb C}$ (recall (\ref{3.11a})). 
Then
\begin{equation}\label{9.2}
{\bf h}=\coprod_{n\in\Z/p{\Bbb Z}}{\bf h}_{(n)}
\end{equation}
(we identify ${\bf h}_{(n\ {\rm mod}\,p)}$ with ${\bf h}_{(n)}$ for
$n\in{\Bbb Z}$). For $\alpha\in{\bf h}$ write $\alpha_{(n)}$ for the
component of $\alpha$ in ${\bf h}_{(n)}.$ 
We define the $\nu$-twisted affine Lie algebra $\hat{\bf h}[\nu]$ 
associated with the abelian Lie algebra ${\bf h}$
to be
\begin{equation}\label{9.3}
\hat{\bf h}[\nu]=\coprod_{n\in\frac{1}{p}{\Bbb
Z}}{\bf h}_{(pn)}\otimes t^{n}\oplus{\Bbb C}c
\end{equation}
with
\begin{equation}\label{9.4}
[x\otimes t^m,y\otimes t^n]=\langle x,y\rangle m\delta_{m+n,0}c\ \ {\rm
for}\ \ x\in{\bf h}_{(pm)},y\in{\bf h}_{(pn)},\ \
m,n\in\frac{1}{p}{\Bbb Z}
\end{equation}
\begin{equation}
[c,\hat{\bf h}[\nu]]=0.
\end{equation}
As in Section 3, this notation and related notation below may be applied
to any finite-dimensional abelian Lie algebra with a nonsingular
symmetric bilinear form and with a finite-order isometry. Note that
for the identity 
automorphism $\nu=1,$ with $p$ chosen to be 1,
the twisted algebra $\hat{\bf h}[\nu]$  reduces
to the untwisted algebra $\hat{\bf h}$ of (\ref{3.1})-(\ref{3.3}). 
Define the {\em weight gradation} on $\hat{\bf h}[\nu]$ associated 
with ${\bf h}_*$ by:
\begin{equation}
{\rm wt}\,(x\otimes t^{m})=0, \ \ {\rm wt}\,(y\otimes t^{n})=-n, \ \ \
{\rm wt}\,c=0
\end{equation}
for $m,n\in \frac{1}{p}{\Bbb Z},$ 
$x\in ({\bf h}_*)_{(pm)},\ y\in ({\bf h}_*^{\perp})_{(pn)}$  
(cf. (\ref{3.5})).  Set 
\begin{equation}
\hat{\bf h}[\nu]^+=\coprod_{n>0}{\bf h}_{(pn)}\otimes t^{n},\ \ 
\hat{\bf h}[\nu]^-=\coprod_{n<0}{\bf h}_{(pn)}\otimes t^{n}.
\end{equation}
Now the subalgebra 
\begin{equation}\label{9.5}
\hat{\bf h}[\nu]_{\frac{1}{p}{\Bbb Z}}=\hat{\bf h}[\nu]^+\oplus
\hat{\bf h}[\nu]^-\oplus {\Bbb C}c
\end{equation}
of $\hat{\bf h}[\nu]$ is a Heisenberg algebra (cf. (\ref{3.7})). Form the 
induced
$\hat{\bf h}[\nu]
$-module 
\begin{equation}\label{9.6}
S[\nu]=U(\hat{\bf h}[\nu])\otimes_{U(\coprod_{n\ge
0}{\bf h}_{(pn)}\otimes t^{n}\oplus{\Bbb C}c)}{\Bbb C}\simeq S(\hat{\bf
h}[\nu]^-)\ \ \ {\rm (linearly)},
\end{equation}
which is irreducible under $\hat{\bf h}[\nu]_{\frac{1}{p}{\Bbb Z}},$ where $\coprod_{n\ge 0}{\bf h}_{(pn)}\otimes t^{n}$ acts
trivially on ${\Bbb C}$ and $c$ acts as 1. We give the module $S[\nu]$ a
${\Bbb Q}$-grading compatible with the action of $\hat{\bf h}[\nu],$
and such that 
\begin{equation}\label{gra}
{\rm wt}\,1=\frac{1}{4p^2}\sum_{k=1}^{p-1}k(p-k){\rm dim}\,({\bf
h}_*^{\perp})_{(k)}
\end{equation}
(the reason for choosing this shifted gradation will be
explained in Section 6).  The lattice-isometry $\nu$
acts naturally on $\hat{\bf h}[\nu]$ (fixing $c$):
$$\nu (\alpha\otimes t^n)=\omega_p^n\alpha\otimes t^n$$ 
for $n\in\frac{1}{p}\Z,\ \alpha\in {\bf h}_{(pn)},$ 
and on $S[\nu]$ as an algebra isomorphism,
preserving the gradation, and we have
\begin{equation}\label{extra4.11}
\nu(u\cdot v)=\nu(u)\cdot\nu(v)
\end{equation}
for $u\in \hat{\bf h}[\nu]$ and $v\in S[\nu].$  For $\alpha \in {\bf h},$ 
$n\in \frac{1}{p}{\Bbb Z}$, write
$\alpha (n)$ for the operator on $S[\nu]$ determined by
$\alpha_{(pn)}\otimes t^n.$ 

Let $T$ be an $\hat L_{\nu}$-module with $\kappa_q$ acting as multiplication 
by $\omega_q.$  We assume that ${\bf h}_{(0)}$ acts on 
$T$ in  
such a way that 
\begin{equation}\label{assume1}
T=\coprod_{\alpha\in{\bf h}_{(0)}\cap{\bf h}_*^{\perp}}T_\alpha
\end{equation}
where
\begin{equation}\label{assume2}
T_\alpha=\{t\in T|h\cdot t=\langle h,\alpha\rangle t\ \ {\rm for\ \ }
h\in {\bf h}_{(0)}\},
\end{equation}
and such that the actions of $\hat L_\nu$ and ${\bf h}_{(0)}$ are 
compatible in the 
sense that
\begin{equation}\label{assume3}
a\cdot T_{\alpha}\subset T_{\alpha+\bar a_{(0)}'}
\end{equation}
for $a\in \hat L_\nu$ and $\alpha\in {\bf h}_{(0)}\cap{\bf h}_*^{\perp}.$
We also assume that $\hat \nu$ acts on $T$ as a linear automorphism such that 
\begin{equation}\label{44.1}
\hat\nu T_\alpha\subset T_{\nu\alpha}.
\end{equation}
Then as operators on $T,$
\begin{equation}\label{assume4}
h\,a=a(\langle h',\bar a\rangle+h)
\end{equation}
\begin{equation}\label{extra4.18}
\hat\nu h\hat\nu^{-1}=h
\end{equation}
for $h\in {\bf h}_{(0)}$ and $a\in\hat L_\nu.$ 
Define a ${\Bbb C}$-gradation on $T$  by
\begin{equation}\label{44.2}
{\rm wt}\, t=\frac{1}{2}\langle \alpha',\alpha'\rangle\ \ \
{\rm for }\ \ \ t\in T_{\alpha'} \ (\alpha\in {\bf h}_{(0)}).
\end{equation}
Then $\hat\nu$ preserves this gradation of $T$ by (\ref{44.1}).
We define an End\,$T$-valued formal Laurent series
$z^h$ for $h\in {\bf h}_{(0)}$ as follows:
\begin{equation}\label{assume5}
z^h\cdot t=z^{\langle h,\alpha\rangle}t \ \ \ {\rm for}\ \ \ t\in
T_{\alpha}\ (\alpha\in{\bf h}_{(0)}\cap{\bf h}_*^{\perp}).
\end{equation}
Then from (\ref{assume4}),
\begin{equation}\label{assume6}
z^ha=az^{\langle h',\bar a\rangle+h}\ \ \ {\rm for}\ \ \ a\in \hat L_{\nu}
\end{equation}
as operators on $T.$ 

\begin{rem}\hspace{-0.2 cm}{\bf : }{\rm If $\nu=1$ and 
$c_0(\cdot,\cdot)=c_0^{\nu}(\cdot,\cdot),$ then $\hat L_{\nu}=\hat L$
and ${\Bbb C}\{L\}$ is a such module $T.$}
\end{rem}

\begin{rem}\label{erem4.4}\hspace*{-0.2 cm}{\bf : }{\rm In the setting of
 Remark \ref{rem1.2},  
let $L$ be an even lattice satisfying (\ref{L2.1}) with the alternating
bilinear map $c_0^{\nu}$ given in (\ref{add1}). Then $K=\{a^{-1}\hat\nu(a)|a\in\hat L_{\nu}\}$ is a central subgroup of $\hat L_{\nu}$ 
and $K\cap \langle \kappa_q\rangle=1.$ 
In [L1], a certain class of 
$\hat L_{\nu}$-modules on which $\kappa_q$ acts by $\omega_q$ and
$K$ acts according to the character $\chi(a^{-1}\hat\nu a)=\omega_p^{-\langle
\sum \nu^r\bar a,\bar a\rangle/2}$ is classified and constructed explicitly
(see Propositions 6.1 and 6.2 of [L1]).  
These modules have the 
properties described above with ${\bf h}_*^{\perp}={\bf h}.$
In particular, any such  
module $T$ (which is denoted by $U_T$ in [L1]) 
has the following decomposition:
\begin{equation}
T=\sum_{a\in L}T_{\bar a_{(0)}}.
\end{equation}
For $\alpha\in{\bf h}_{(0)},$ if $\langle \alpha,\bar a_{(0)}\rangle\in {\Bbb Z}$  
for all $a\in L,$ define the operator $\omega_q^{\alpha}$ on $T$ by
\begin{equation}\label{oct1}
\omega_q^{\alpha}\cdot t=\omega_q^{\langle\alpha,\bar a_{(0)}\rangle}t
\end{equation}
for $t\in T_{\bar a_{(0)}}.$ Then as operators on $T,$
\begin{equation}\label{right}
\hat\nu a=a\omega_p^{-p\bar a_{(0)}-p\langle\bar a_{(0)},\bar a_{(0)}\rangle/2},
\end{equation}
where $a\in \hat L_{\nu}$ and the operator in the right-hand side  is
well defined because of the assumption (\ref{L2.1}). It follows that
\begin{equation}
\hat\nu^p=1
\end{equation}
(see [L1]), a nontrivial fact that is not an automatic consequence of the assumption $\nu^p=1.$}
\end{rem}
 
Set
\begin{equation}\label{9.8}
V^T_L=S[\nu]\otimes T
\end{equation}
which is naturally graded, using the gradations of $S[\nu]$ and $T.$ 
Again  $\hat{L}_{\nu},$
$\hat{\bf h}[\nu]_{\frac{1}{p}{\Bbb Z}},$ ${\bf h}_{(0)},$ $z^h$
$(h\in{\bf h}_{(0)})$ act on $V_L^T$ by acting on either $S[\nu]$ or
$T$. Then $\hat\nu$ extends to a linear automorphism of $V_L^T$ so that
$\hat\nu(u\otimes t)=\nu(u)\otimes \hat\nu(t)$ for $u\in S[\nu]$ and $t\in
T.$ As in Section 3, we write
\begin{equation}
c_{\nu}(\alpha,\beta)=\omega_q^{c_0^{\nu}(\alpha,\beta)}
\end{equation}
for $\alpha,\beta\in L.$ Then for $a,b\in \hat{L}_{\nu}$  
\begin{equation}
aba^{-1}b^{-1}=c_{\nu}(\bar a,\bar b)
\end{equation}
as operators on $T$ and on $V_L^T.$ 
We define two subspaces of $V_L^T:$  
\begin{equation}\label{9.9}
\Omega^{\nu}_*=\{v\in V^T_L\vert h(n)v=0,\ \mbox{for}\ \ h\in {\bf h}_*,
n\in\frac{1}{p}{\Bbb Z}, n>0\}
\end{equation}
\begin{equation}\label{9.10}
V^{\nu}_*=\{h(n)V^T_L\vert h\in {\bf h}_*, n\in\frac{1}{p}{\Bbb Z},
n<0\}
\end{equation}
(cf. (\ref{3.29}) and (\ref{3.30})). Then $\Omega^{\nu}_*$ is the vacuum space for the Heisenberg algebra 
${\bf h}_*[\nu]_{\frac{1}{p}{\Bbb Z}}$ and
\begin{equation}\label{9.11}
V^T_L=V^{\nu}_*\oplus \Omega^{\nu}_*.
\end{equation}
As in (\ref{3.33}) and (\ref{3.34}) we have
\begin{equation}
\Omega_*^{\nu}=S((\hat{\bf h}_*^{\perp}[\nu])^-)\otimes T
\end{equation}
\begin{equation}
V_*^{\nu}=\hat{\bf h}_*[\nu]^-S(\hat{\bf h}_*[\nu]^-)\otimes\Omega_*^{\nu}.
\end{equation}
The linear map $\hat\nu$ preserves both $\Omega_*^{\nu}$ and $V_*.$

For a nonzero complex number $a$ we shall fix the branch of
$a^{\tau}=e^{\tau(\log\vert a\vert +i\,{\rm arg}\,a)}$ ($\tau$ is a complex
variable) with $-\pi<\,$arg\,$a\leq\pi.$ Then $a^{\tau+\zeta}=
a^{\tau}a^{\zeta}$ for $\tau,\zeta\in {\Bbb C}.$ The following lemma will
be used in proving the Jacobi identity in the next section:
\begin{lemma}\label{ltau} For $\tau\in{\Bbb C},$ 
\begin{equation}\label{tau}
\prod_{r=1}^{p-1}(1-\omega_p^r)^{\tau}=p^{\tau}.
\end{equation}
\end{lemma}

\noindent{\bf Proof\ } If $\omega_p^r\ne -1,$ arg$(1-\omega_p^r)=-$arg$(1-\omega_p^{-r}).$ Then
$$(1-\omega_p^r)^{\tau}(1-\omega_p^{-r})^{\tau}=((1-\omega_p^r)(1-\omega_p^{-r}))^{\tau}$$
and the lemma follows from the fact that $\prod_{r=1}^{p-1}(1-\omega_p^r)=p.$
\ \ $\Box$

For $\alpha\in L$, we define
\begin{equation}\label{9.12}
\sigma(\alpha)=\left\{\begin{array}{ll}
{\displaystyle \prod_{0<r<p/2}(1-\omega^{-r}_p)^{\langle\nu ^r\alpha',
\alpha'\rangle}
2^{\langle\nu^{p/2}\alpha',\alpha'\rangle}} & {\rm if}\ \
p\in 2{\Bbb Z}\\
{\displaystyle\prod_{0<r<p/2}(1-\omega^{-r}_p)^{\langle\nu ^r\alpha',\alpha'\rangle}}
& {\rm if}\ \ p\in 2{\Bbb Z}+1.\end{array}\right.
\end{equation}
Then $\sigma(\nu\alpha)=\sigma(\alpha)$. 
We also define
\begin{equation}
\epsilon_1(\alpha,\beta)=\prod_{0<r<p}(1-\omega_p^{-r})^{\langle\nu^r\alpha,
\beta\rangle}
\end{equation}
\begin{equation}
\epsilon_2(\alpha,\beta)=\omega_q^{\epsilon_0(\alpha,\beta)}
\end{equation}
for $\alpha,\beta\in L.$  

\begin{rem}\hspace*{-0.2 cm}{\bf :\ }{\rm  
In the situation of Remark 2.2, the 
functions $\sigma(\alpha),$ $\epsilon_1(\alpha,\beta)$
and $\epsilon_2(\alpha,\beta)$ 
were introduced in [L1] with ${\bf h}_*^{\perp}={\bf h}$ to define general twisted vertex 
operators. One can check that in this case
\begin{equation}\label{ep}
\epsilon_1(\alpha,\beta)=\epsilon_2(\alpha,\beta)\sigma(\alpha+\beta)/
\sigma(\alpha)\sigma(\beta)
\end{equation}
for $\alpha,\beta\in L.$  }
\end{rem}

For $\alpha\in {\bf h}$,
set
\begin{equation}\label{9.13}
\alpha(z)=\sum_{n\in\frac{1}{p}{\Bbb Z}}\alpha(n)z^{-n-1}.
\end{equation}
We now define the $relative$ $\nu$-$twisted$ $vertex$ $operator$ 
$Y^{\nu}_*(a,z)$ for $a\in\hat{L}$ acting on $V_L^T$ as follows:
\begin{equation}\label{9.14}
Y^{\nu}_*(a,z)=p^{-\langle \bar{a}',\bar{a}'\rangle /2}\sigma(\bar{a})\mbox{$\circ\atop\circ$}
e^{\int(\bar{a}'(z)-\bar{a}'(0)z^{-1})}az^{\bar{a}'_{(0)}+\langle \bar{a}'_{(0)}
,\bar{a}'_{(0)}\rangle /2-\langle \bar{a}',\bar{a}'\rangle /2}  \mbox{$\circ\atop\circ$}
\end{equation}
generalizing the relative untwisted vertex operator $Y_*(a,z)$, the
case $\nu =1$. Note that on the right-hand side, we view $a$ as an
element of $\hat L_{\nu},$ according to 
our set-theoretic identification between
$\hat L$ and $\hat L_{\nu}$ (see Remark \ref{rem1.1}). The numerical factor
at the front of (\ref{9.14}) leads to just the right form for the
general result below.

For $\alpha_1,...,\alpha_k \in{\bf h},$ $n_1,...,n_k
\in {{\Bbb Z}}\ (n_i>0)$ and
$v=\alpha_1(-n_1)\cdot\cdot\cdot\alpha_k(-n_k)\cdot\iota(a)\in V_L,$
set
\begin{equation}\label{9.15}
W_*(v,z)=\mbox{$\circ\atop\circ$}\left(\frac{1}{(n_1-1)!}\left(\frac{d}{dz}\right)^{n_1-1}\alpha'_1(z)\right)\cdot\cdot\cdot\left(\frac{1}{(n_k-1)!}\left(\frac{d}{dz}\right)^{n_k-1}
\alpha'_k(z)\right)Y^{\nu}_*(a,z)\mbox{$\circ\atop\circ$},
\end{equation}
where the right side is an operator on $V^T_L$. Extend to all $v\in V_L$ by
linearity. 

Let $\{\beta_1,\cdot\cdot
\cdot, \beta_d\}$ be an orthonormal basis of
${\bf h}^{\perp}_*,$  and define constants $c_{mni}\in{\Bbb C}$
for $m, n\ge0$ and $i=0,\cdot\cdot\cdot, p-1$ by the formulas
\begin{equation}\label{cmn}
\begin{array}{c}
{\displaystyle\sum_{m,n\ge 0}c_{mn0}x^my^n=-\frac{1}{2}\sum_{r=1}^{p-1}{\rm log}\left(\frac
{(1+x)^{1/p}-\omega^{-r}_p(1+y)^{1/p}}{1-\omega^{-r}_p}\right),}\\
{\displaystyle\sum_{m,n\ge 0}c_{mni}x^my^n=\frac{1}{2}{\rm log}\left( \frac
{(1+x)^{1/p}-\omega^{-i}_p(1+y)^{1/p}}{1-\omega^{-i}_p}\right)}\ \
{\rm for}\ \ i\ne0.
\end{array}
\end{equation}
Set
\begin{equation}\label{9.16}
\Delta_z=\sum_{m,n\ge 0}\displaystyle{\sum^{p-1}_{i=0}}\ \displaystyle{
\sum^d_{j=1}}
c_{mni}(\nu^{-i}\beta_j)(m)\beta_j(n)z^{-m-n}.
\end{equation}
Then $e^{\Delta_z}$ is well defined on $V_L$ since\ \ $c_{00i}=0$ for all
$i$, and for $v\in V_L,$ $e^{\Delta_z}v\in V_L[z^{-1}].$ Note that
$\Delta_z$ is independent of the choice of orthonormal basis. 
Then from (\ref{e3.48}),
\begin{equation}
\hat\nu\Delta_z=\Delta_z\hat\nu
\end{equation}
and hence
\begin{equation}\label{44.4}
\hat\nu e^{\Delta_z}=e^{\Delta_z}\hat\nu
\end{equation}
on $V_L.$ For $v\in V_L,$ the $relative$ $\nu$-$twisted$ $vertex$ $operator$
$Y^{\nu}_*(v,z)$ is defined by:
\begin{equation}\label{9.17}
Y^{\nu}_*(v,z)=W_*(e^{\Delta_z}v,z).
\end{equation}
Then this yields a well-defined linear map
\begin{equation}\label{9.18}
\begin{array}{lcr}
V_L &\to&(\mbox{End}\,V^T_L)\{z\}\hspace*{4.8 cm}\\
\ v &\mapsto& Y^{\nu}_*(v,z)=\displaystyle{\sum_{n\in{\Bbb
C}}v^{\nu}_nz^{-n-1}\ \ (v^{\nu}_n\in {\rm End}\,V^T_L).}
\end{array}
\end{equation}

We now discuss the relation between $\nu$ and twisted vertex operators
$Y_*^{\nu}(v,z)$ for $v\in V_L.$ It is easy to see from the
definitions (\ref{9.14}) and (\ref{9.15}) that
\begin{equation}
\nu\alpha(n)\nu^{-1}=(\nu\alpha)(n)\ \ \ {\rm for}\ \ \ \alpha\in{\bf
h},\ n\in\frac{1}{p}{\Bbb Z}
\end{equation}
\begin{equation}
\nu Y_*^{\nu}(a,z)\nu^{-1}=Y_*^{\nu}(\hat{\nu}a,z)\ \ \ {\rm for}\ \ \
a\in \hat L
\end{equation}
\begin{equation}
\nu W_*(v,z)\nu^{-1}=W_*(\hat\nu v,z)\ \ \ {\rm for}\ \ \ v\in V_L
\end{equation}
and it follows from (\ref{44.4}) and (\ref{9.17}) that
\begin{equation}
\nu Y_*^{\nu}(v,z)\nu^{-1}=Y_*^{\nu}(\hat{\nu}v,z)\ \ \ {\rm for}\ \ \
v\in V_L.
\end{equation}
In particular, $\nu$ commutes with  $Y_*^{\nu}(v,z)$ if $v$ is $\hat
\nu$-invariant. By the properties of $\hat\nu$ and its action on $V_L^T,$
we also have:
\begin{equation}\label{extral}
Y_*^{\nu}(\hat\nu^rv,z)=\lim_{z^{1/p}\to\omega_p^{-r}z^{1/p}}Y_*^{\nu}(v,z)
\omega_p^{pr\bar a'_{(0)}+pr\langle\bar a'_{(0)},\bar a'_{(0)}\rangle/2-
pr\langle\bar a',\bar a'\rangle/2}a^{-1}\hat\nu^ra
\end{equation}
for $v=v'\otimes\iota(a)\in V_L$ where $v'\in S(\hat{\bf h}^-),$ 
$a\in\hat L$ and the operator $\omega_p^{\alpha}$ is defined as 
in (\ref{oct1}). Note that $\omega_p^{\tau}$ is well defined for any
complex number $\tau$.

\begin{rem}\label{erem4.2}\hspace*{-0.2 cm}{\bf : }{\rm 
In the context of Remark 2.2, 
let $T$ be an $\hat L_{\nu}$-module as
given in Remark \ref{erem4.4}. Then
the operators $Y^{\nu}_*(v,z)$ are exactly
the $\nu$-twisted vertex operators $Y_{\nu}(v,z)$ introduced in [FLM2]. In
particular, $Y_{\nu}(u,z)\in ({\rm End}\,V_L^T)[[z^{1/p},z^{-1/p}]]$ for
$u\in V_L$ and (\ref{extral}) reduces to
\begin{equation}\label{e4.53}
Y_{\nu}(\hat\nu^rv,z)=\lim_{z^{1/p}\to\omega_p^{-r}z^{1/p}}Y_{\nu}(v,z)
\end{equation}
(see (\ref{right})).}
\end{rem}

We summarize the main elementary properties of relative twisted vertex 
operators in
the following proposition.
\begin{prop}
Let $a, b\in\hat{L},$ $v\in V_L,$ $\alpha\in{\bf h}_{(0)},$
$\psi\in\hat{L}$ such that $\bar{\psi}\in L\cap {\bf h}_*$. We have
\begin{equation}\label{9.19}
\mbox{$\circ\atop\circ$}Y^{\nu}_*(a,z_1)Y^{\nu}_*(b,z_2)
\mbox{$\circ\atop\circ$}
=c_{\nu}(a,b)\mbox{$\circ\atop\circ$}Y^{\nu}_*(b,z_2)Y^{\nu}_*(a,z_1)
\mbox{$\circ\atop\circ$}
\end{equation}
\begin{equation}\label{9.20}
Y^{\nu}_*(a,z_1)Y^{\nu}_*(b,z_2)=\mbox{$\circ\atop\circ$}
Y^{\nu}_*(a,z_1)Y^{\nu}_*(b,z_2)\mbox{$\circ\atop\circ$}
\displaystyle{\prod^{p-1}_{i=0}}(z^{1/p}_1-\omega^{-i}_pz^{1/p}_2)^{\langle \nu^i\bar{a}',\bar{b}'\rangle }
\end{equation}
\begin{equation}\label{9.21}
\frac{d}{dz}Y^{\nu}_*(a,z)=\mbox{$\circ\atop\circ$}
\left(\bar{a}'(z)+\frac{1}{2}\langle \bar{a}'_{(0)},\bar{a}'_{(0)}\rangle z^{-1}-\frac{1}{2}
\langle \bar{a},\bar{a}\rangle z^{-1}\right)Y^{\nu}_*(a,z)\mbox{$\circ\atop\circ$}
=Y^{\nu}_*(\bar{a}(-1)a,z)
\end{equation}
\begin{equation}\label{9.22}
Y^{\nu}_*(v,z)=0\ \ \mbox{if} \ \ v\in V_*
\end{equation}
\begin{equation}\label{9.23?r}
[\hat{\bf h}_*[\nu]_{\frac{1}{p}{\Bbb Z}},Y^{\nu}_*(v,z)]=0
\end{equation}
and $v^{\nu}_n$ preserves both $\Omega^{\nu}_*$ and
$V^{\nu}_*$ for $n\in {\Bbb C},$ 
\begin{equation}\label{9.24} 
Y^{\nu}_*(\psi v,z)=\omega_q^{\epsilon_0(\bar\psi,\bar a)}\psi Y^{\nu}_*(v,z)
\end{equation}
if $v=v^*\otimes \iota(a)$ for $v^*\in S(\hat{\bf h}^-)$ (see (\ref{2.6})),
and in particular
\begin{equation}\label{9.25}
Y^{\nu}_*(\iota({\psi}),z)=\psi.
\end{equation}
(The binomial expressions in (\ref{9.20})  
are to be expanded in nonnegative
integral powers of the second variable $z^{1/p}_2;$ recall Remark 3.3.)
\end{prop}

\noindent{\bf Proof } First, (\ref{9.19}) follows from the facts that
$$\mbox{$\circ\atop\circ$}Y^{\nu}_*(a,z_1)Y^{\nu}_*(b,z_2)
\mbox{$\circ\atop\circ$}=p^{-\langle \bar{a}',\bar{a}'\rangle /2-\langle \bar{b}',\bar{b}'\rangle /2}
\sigma(\bar{a})\sigma(\bar{b})
\mbox{$\circ\atop\circ$}e^{\int(\bar{a}'(z_1)-\bar{a}'(0)z^{-1}_1
+\bar{b}'(z_2)-\bar{b}'(0)z^{-1}_2)}\mbox{$\circ\atop\circ$}\cdot$$
\begin{equation}\label{9.26}
\cdot abz_1^{\bar{a}'_{(0)}+\langle \bar{a}'_{(0)},\bar{a}'_{(0)}\rangle /2-
\langle \bar{a}',\bar{a}'\rangle /2}  
z_2^{\bar{b}'_{(0)}+\langle \bar{b}'_{(0)},\bar{b}'_{(0)}\rangle /2-
\langle \bar{b}',\bar{b}'\rangle /2}  
\end{equation}
and that $ab=c_{\nu}(a,b)ba.$ (\ref{9.20}) follows from (\ref{assume6})
and
$$e^{\int(\bar{a}'(z_1)-\bar{a}'(0)z^{-1}_1
+\bar{b}'(z_2)-\bar{b}'(0)z^{-1}_2)}=
\mbox{$\circ\atop\circ$}e^{\int(\bar{a}'(z_1)-\bar{a}'(0)z^{-1}_1
+\bar{b}'(z_2)-\bar{b}'(0)z^{-1}_2)}\mbox{$\circ\atop\circ$}\cdot $$
$$\cdot
\displaystyle{\prod^{p-1}_{i=0}}\left( 1-\omega^{-i}_p\left( \frac{z_2}{z_1}
\right)^{1/p}\right)
^{\langle \nu^i\bar{a}',\bar{b}'\rangle },$$
and (\ref{9.21}) follows from the computation
$$e^{\Delta_z}\bar{a}(-1)a=
\left(\bar a'(-1)+\frac{1}{2}\langle \bar{a}'_{(0)},\bar{a}'_{(0)}
\rangle z^{-1}-\frac{1}{2}\langle \bar{a}',\bar{a}'\rangle z^{-1}\right)a$$
and the definition of relative twisted vertex operators (\ref{9.17}). The other
relations are the direct consequences of the definitions (\ref{9.14}),
(\ref{9.15}) and (\ref{9.17}). \ \ \ \ \ \ \ $\Box$

Since $Y_*^{\nu}(u,z)=0$ for $u\in V_*$ we shall be especially interested in
the operators $Y_*^{\nu}(v,z)$ for $v\in\Omega_*.$ We also use the notation $Y_*^{\nu}$ for the restriction map:
\begin{equation}
\begin{array}{lll}
\Omega_* &\to&(\mbox{End}\,V^T_L)\{z\}\\
\ v &\mapsto& Y^{\nu}_*(v,z)=\displaystyle{\sum_{n\in{\Bbb
C}}v^{\nu}_nz^{-n-1}.}
\end{array}
\end{equation}

\section{A Jacobi identity for the relative twisted vertex
operators}

\setcounter{equation}{0}

We continue our discussion of relative vertex operators. Here we 
present the main theorem $-$ a Jacobi identity for relative
twisted vertex operators, whose proof closely follows the pattern of
the proof of the corresponding result in Chapter 9 of [FLM3]. 
This result generalizes a large number of known ones.
The case $\nu=1$ amounts to the theory of relative untwisted vertex
operators of [DL1]-[DL2] in which we clarified the essential
equivalence
between untwisted 
$Z$-algebras [LP2] and parafermion algebras [ZF1].
If the lattice is even and ${\bf h}_*=0,$ we recover
the twisted vertex operators and the associated results obtained in
[L1] and [FLM2]. In certain special cases we also recover the $Z$-algebra
structure and relations found in [LW1]-[LW5] and [LP1]-[LP3]. 
In the case that $L$ is a direct sum of several
copies of the root lattice of the Lie algebra $sl(2,{\Bbb C})$ and 
$\nu=-1,$ the Jacobi identity for certain relative twisted vertex 
operators was established in [Hu] to explain   
the essential equivalence
between the ``twisted $Z$-algebras'' associated with the principal 
representations of $A_1^{(1)}$ [LW3]-[LW4] and
the ``twisted parafermion algebras'' of [ZF2].

We shall use the basic generating function
\begin{equation}
\delta(z)=\sum_{n\in{\Bbb Z}}z^n,
\end{equation}
formally the expansion of the $\delta$-function at $z=1.$ The fundamental
(and elementary) properties of the $\delta$-function can be found in [FLM3].
The following result concerning formal calculus and the $``\delta$-function'' will be used in the proof of the Jacobi identity:
\begin{prop}\label{10.1}
Let $V$ be a vector space, $a(m,n)$ $(m,n\in\frac{1}{p}{\Bbb Z})$ a
family of operators on $V.$ Let $s,t\in {\Bbb C}$ and let
$\zeta$ be a $p^{th}$ root of unity. Consider the formal
series
\begin{equation}
A(z_1^{1/p},z_2^{1/p})=z_1^sz_2^t\left(\sum_{m,n\in \frac{1}{p}{\Bbb Z}}a(m,n)z_1^mz_2^n\right)
\end{equation}
and assume that
\begin{equation}
\lim_{z_1^{1/p}\to z_2^{1/p}}A(z_1^{1/p},z_2^{1/p}) \ \ {\rm exists},
\end{equation}
that is, for every $k\in {\Bbb Q}$ and $v\in V,$ $a(m,k-m)v=0$ for all but
a finite number of $m\in\frac{1}{p}{\Bbb Z}.$  
Then
\begin{eqnarray}
& &\displaystyle{A(z_1^{1/p},z_2^{1/p})e^{z_0\frac{\partial}{\partial z_1}}
\delta\left(\zeta\left(\frac{z_1}{z_2}\right)
^{1/p}\right)}\nonumber\\
& &\ \ \ =\displaystyle{\zeta^{sp}A(\zeta^{-1}(z_2-z_0)^{1/p},z_2^{1/p})e^{z_0\frac{\partial}{\partial z_1}}
\left( \left(\frac{z_1}{z_2}\right)^s\delta\left(\zeta\left(\frac{z_1}{z_2}\right)^{1/p}\right)\right)}\\
& &\ \ \ =(z_2-z_0)^sz_2^t\left(\sum_{m,n\in \frac{1}{p}{\Bbb Z}}a(m,n)(z_2-z_0)^mz_2^n\right)e^{z_0\frac{\partial}{\partial z_1}}
\left( \left(\frac{z_1}{z_2}\right)^s\delta\left(\zeta\left(\frac{z_1}{z_2}\right)^{1/p}\right)\right).\nonumber
\end{eqnarray}
\end{prop}
\noindent{\bf Proof } We have: 
\begin{eqnarray*}
& &A(z_1^{1/p},z_2^{1/p})e^{z_0\frac{\partial}{\partial z_1}}
\delta\left(\zeta\left(\frac{z_1}{z_2}\right)^{1/p}\right) \\
& &\ \ \ \ \ =\left(\frac{z_1}{z_2}\right)^s
\left(\frac{z_1}{z_2}\right)^{-s}A(z_1^{1/p},z_2^{1/p})
e^{z_0\frac{\partial}{\partial z_1}}\delta\left(\zeta\left(\frac{z_1}{z_2}\right)^{1/p}\right) \\
& &\ \ \ \ \ =\left(\frac{z_1}{z_2}\right)^se^{z_0\frac{\partial}{\partial z_1}}
\left(\left(\frac{z_1-z_0}{z_2}
\right)^{-s}A((z_1-z_0)^{1/p},z_2^{1/p})\delta\left(\zeta\left(\frac{z_1}{z_2}\right)^{1/p}\right)\right)
\\
& &\ \ \ \ \ \ =\left(\frac{z_1}{z_2}\right)^se^{z_0\frac{\partial}{\partial z_1}}
\left(\zeta^{ps}\left(\frac{z_2-z_0}{z_2}
\right)^{-s}A(\zeta^{-1}(z_2-z_0)^{1/p},z_2^{1/p})
\delta\left(\zeta\left(\frac{z_1}{z_2}\right)^{1/p}\right)\right) \\
& &\ \ \ \ \ \ =\zeta^{ps}A(\zeta^{-1}(z_2-z_0)^{1/p},z_2^{1/p})
\left(\frac{z_1}{z_2}\right)^se^{z_0\frac{\partial}{\partial z_1}}
\left( \left(1-\frac{z_0}{z_2}
\right)^{-s}\delta\left(\zeta\left(\frac{z_1}{z_2}\right)^{1/p}\right)\right) \\
& &\ \ \ \ \ \ =\zeta^{ps}A(\zeta^{-1}(z_2-z_0)^{1/p},z_2^{1/p})
\left(\frac{z_1}{z_2}\right)^se^{z_0\frac{\partial}{\partial z_1}}
\left( \left(1-\frac{z_0}{z_1}
\right)^{-s}\delta\left(\zeta\left(\frac{z_1}{z_2}\right)^{1/p}\right)\right) \\
& &\ \ \ \ \ \ =\zeta^{ps}A(\zeta^{-1}(z_2-z_0)^{1/p},z_2^{1/p})
\left(\frac{z_1}{z_2}\right)^se^{z_0\frac{\partial}{\partial z_1}}
\left( \left(\frac{z_1}{z_1-z_0}
\right)^{s}\delta\left(\zeta\left(\frac{z_1}{z_2}\right)^{1/p}\right)\right) \\
& &\ \ \ \ \ \ =\zeta^{ps}A(\zeta^{-1}(z_2-z_0)^{1/p},z_2^{1/p})
\left(\frac{z_1}{z_2}\right)^s\left(\frac{z_1+z_0}{z_1}\right)^se^{z_0\frac{\partial}{\partial z_1}}
\delta\left(\zeta\left(\frac{z_1}{z_2}\right)^{1/p}\right) \\
& &\ \ \ \ \ \ =\zeta^{ps}A(\zeta^{-1}(z_2-z_0)^{1/p},z_2^{1/p})
e^{z_0\frac{\partial}{\partial z_1}}
\left( \left(\frac{z_1}{z_2}\right)^s\delta\left(\zeta\left(\frac{z_1}{z_2}\right)^{1/p}
\right)\right),
\end{eqnarray*}
as desired. \ \ \ \ \ \ $\Box$

For notational convenience we write
\begin{equation}\label{nota}
F_{(\alpha,\beta)}(z_1,z_2)=\prod_{r=0}^{p-1}(z_1^{1/p}-\omega^{-r}_pz_2^{1/p}
)^{-\langle \nu^r\alpha',\beta'\rangle }
\end{equation}
\begin{equation}\label{gg}
G_{(\alpha,\beta)}(z_1,z_2)=\prod_{r=0}^{p-1}(z_1^{1/p}-\omega^{-r}_pz_2^{1/p})^{\langle\alpha',\beta'\rangle }F_{(\alpha,\beta)}(z_1,z_2)
\end{equation}
for $\alpha,\beta\in{L}$, where as usual all binomial expressions are to be 
expanded in nonnegative 
integral powers of the second variable. We shall often use 
expressions like $F_{\a,\b}(z_2+z_0,z_2)$ and $G_{\a,\b}(z_2+z_0,z_2)$
below. It is understood that the expression 
$(z_2+z_0)^{1/p}-\omega_p^rz_2^{1/p}$ 
is to be expanded in nonnegative integral powers of $z_0,$ and so we
have expansions of the form: 
$$
((z_2\!+\!z_0)^{1/p}\!-\!\omega^{r}_pz_2^{1/p})^c\!
=\!\left\{\begin{array}{ll}
p^{-c}z_2^{c/p}(z_0/z_2)^c\left(\!1\!+\!\displaystyle{\sum_{n>0}}a_n(z_0/z_2)^n\right)
 & {\rm if}\ r=0,\\
 z_2^{c/p}(1\!-\!\omega_p^r)^{c}\left(\!1\!+\!\displaystyle{\sum_{n>0}}b_n
(z_0/z_2)^n\right) & {\rm if}\ 
r>0.
\end{array}\right.$$
We also write
\begin{equation}\label{e5.8}
\tau(\alpha,\beta)=\frac{\sigma(\alpha)\sigma(\beta)\epsilon_1(\alpha,\beta)}
{\sigma(\alpha+\beta)\epsilon_2(\alpha,\beta)}\in {\Bbb C}
\end{equation}
for $\alpha,\beta\in L.$ 
Note that in the context of Remark 2.2, 
$\alpha'=\alpha,$ $\beta'=\beta$ and  
$\tau(\alpha,\beta)=1$, by (\ref{ep}). The following general residue notation
will be useful:
\begin{equation}
{\rm Res}_z(\sum v_nz^n)=v_{-1}.
\end{equation}

\begin{th}\label{t10.2}
Let $a, b\in\hat{L},$ $u^*, v^*\in M(1)$
and set
\begin{equation}
\begin{array}{l}
u=u^*\otimes \iota(a)\in V_L\\
v=v^*\otimes\iota(b)\in V_L.
\end{array}
\end{equation}
Then we have:
$$F_{(\bar a,\bar b)}(z_1,z_2)(z_1-z_2)^nY^{\nu}_*(u,z_1)Y^{\nu}_*(v,z_2)
-c_{\nu}(\bar{a},\bar{b})F_{(\bar b,\bar a)}(z_2,z_1)(-z_2+z_1)^n
Y^{\nu}_*(v,z_2)Y^{\nu}_*(u,z_1)$$
\begin{eqnarray}
&={\rm Res}_{z_0}\frac{1}{p}z_2^{-1}z^n_0{\displaystyle\sum^{p-1}_{r=0}}
\left(G_{(\nu^r\bar a,\bar b)}(z_2+z_0,z_2)z_0^{-\langle\nu^r\bar a',\bar b'\rangle}
Y^{\nu}_*(Y_*(\hat{\nu}^ru,z_0)v,z_2)\tau(\nu^r\bar a,\bar
b)\cdot\right.\\
&\ \ \ \cdot\left.c_{\nu}(\bar a-\nu^r\bar a,\bar b)(\hat\nu^ra^{-1})a
e^{-z_0\frac{\partial}{\partial z_1}}
\delta(\omega^r_pz_1^{1/p}/z_2^{1/p})
(z_1/z_2)^{\bar{a}'_{(0)}+\langle \bar{a}'_{(0)},\bar{a}'_{(0)}\rangle /2
-\langle \bar{a}',\bar{a}'\rangle /2}\right)\nonumber 
\end{eqnarray}
for $n\in{\Bbb Z},$ or equivalently,
$$z^{-1}_0F_{(\bar a,\bar b)}(z_1,z_2)\delta\left(\frac{z_1-z_2}{z_0}\right)
Y^{\nu}_*(u,z_1)Y^{\nu}_*(v,z_2)$$
\begin{equation}\label{tjac}
-c_{\nu}(\bar{a},\bar{b})z^{-1}_0
F_{(\bar b,\bar a)}(z_2,z_1)
\delta\left(\frac{z_2-z_1}{-z_0}\right)
Y^{\nu}_*(v,z_2)Y^{\nu}_*(u,z_1)
\end{equation}
$$=\frac{1}{p}z_2^{-1}\displaystyle{\sum^{p-1}_{r=0}}
\left(G_{(\nu^r\bar a,\bar b)}(z_2+z_0,z_2)
z_0^{-\langle\nu^r\bar a',\bar
b'\rangle}Y^{\nu}_*(Y_*(\hat{\nu}^ru,z_0)v,z_2)\delta\left(\omega^r_p\left(\frac{z_1-z_0}{z_2}\right)^{1/p}\right)\cdot\right.
$$
$$\left.\cdot\tau(\hat\nu^r\bar a,\bar b)c_{\nu}(\bar a-\nu^r\bar a,\bar b)(\hat\nu^{r}a^{-1})a\left(\frac{z_1-z_0}{z_2}\right)^{\bar{a}'_{(0)}+\langle \bar{a}'_{(0)},\bar{a}'_{(0)}\rangle /2
-\langle \bar{a}',\bar{a}'\rangle /2}\right),$$
where all binomial expressions
are to be expanded in nonnegative integral powers of the second
variables.
\end{th}

\noindent {\bf Proof } Let $k, l\ge 1$ and let
$a_1,\cdot\cdot\cdot,a_k, b_1,\cdot\cdot\cdot,b_l\in\hat{L}$ subject
to the conditions $a=a_1\cdots a_k$ and $b=b_1\cdots b_l.$
Then the coefficients in the formal power series $A$ and
$B$
\begin{equation}\label{main1}
\begin{array}{c}
{\displaystyle 
A={\rm exp}\left(\sum_{i=1}^{k}\sum_{n>0}\frac{\bar a_i'(-n)
}{n}w_i^n\right)\iota(a)\in V_L[[w_1,...,w_k]]}\\
{\displaystyle 
B={\rm exp}\left(\sum_{j=1}^{l}\sum_{n>0}\frac{\bar b_j'(-n)
}{n}x_j^n\right)\iota(b)\in V_L[[x_1,...,x_l]]}
\end{array}
\end{equation}
span $S(({\bf h}^{\perp}_*)^-)\otimes\iota(a)$
and $S(({\bf h}^{\perp}_*)^-)\otimes\iota(b),$
respectively, and so it suffices to prove the theorem 
with $u$ and $v$ replaced by $A$
and $B,$ respectively. 

Let
$\alpha_1,...,\alpha_s
\in L$ and  $\alpha=\alpha_1+\cdot\cdot\cdot+\alpha_s.$ For brevity, set
\begin{equation}\label{5.21}
f_\alpha(z,w_1,...,w_s)=\prod_{
1\le i<j\le s}\,
\prod_{0\le r<p}\left
((z+w_i)^{1/p}-\omega^{-r}_p(z+w_j)^{1/p}\right
)^{\langle (\nu^r-1)\alpha'_i,\alpha'_j\rangle}
\end{equation}
\begin{equation}\label{5.22}
\begin{array}{l}
f_{(\bar a,\bar b)}(z_1,w_1,...,w_k;z_2,x_1,...,x_l)\\
\displaystyle{\ \ \ \ \ \ =\prod_{1\le i\le k,\, 1\le j\le l}\,
\prod_{0\le r<p}}
\left((z_1+w_i)^{1/p}-\omega^{-r}_p(z_2+x_j)^{1/p}\
\right)^{\langle \nu^r\bar{a}'_i,\bar{b}'_j\rangle },
\end{array}
\end{equation}
where all binomial expressions are to be understood as formal power
series in the $w$'s and $x$'s, and the coefficient of each monomial in the $w$'s 
in $f_{\alpha}$ is the product of 
$z^{\langle \bar a'_{(0)}-\bar a',\bar b\rangle}$ and a Laurent polynomial
in $z,$ and the coefficient of each monomial in the $w$'s and $x$'s in $f_{\alpha}$ 
is the product of $z_1^{\langle \bar a_{(0)}',\bar b\rangle}$ 
with a Laurent polynomial in $z_1^{1/p}$ and $z_2^{1/p}.$ We shall also use
$f_{(\bar b,\bar a)}(z_2,x_1,...,x_l;z_1,w_1,...,w_k).$ 

By (\ref{2.6}), we see that as operators on $V_L^T$
\begin{equation}
a=\prod_{1\leq i<j\leq k}\epsilon_2(\bar a_i,\bar a_j)a_1\cdots a_k.
\end{equation}
Just as in the proof of Theorem 9.3.1 of [FLM3], we have 
\begin{eqnarray}
& &W_*(A,z)=\mbox{$\circ\atop\circ$}{\displaystyle{\rm exp}\left(\sum_{i=1}^{k}\sum
_{n\ge1}\frac{1}{n!}\left(\frac{d}{dz}\right)^{n-1}\bar{a}'_i(z)w^n_i\right)Y^{\nu}_*(a,z)} 
\mbox{$\circ\atop\circ$}\hspace{5 cm} \nonumber\\
& &\ \ \ \ =\mbox{$\circ\atop\circ$}Y^{\nu}_*(a_1,z+w_1)\cdot\cdot\cdot
Y^{\nu}_*(a_k,z+w_k)\mbox{$\circ\atop\circ$}\prod_{1\le i<j\le k}p^{-\langle \bar{a}_i',\bar{a}_j'\rangle}\epsilon_2(\bar a_i,\bar a_j)\cdot\label{5.23}\\
& &\ \ \ \ \cdot\frac{\sigma(\bar a)}{\sigma(\bar
a_1)\cdots\sigma(\bar a_k)}\prod_{1\le i<j\le
k}z^{\langle (\bar{a}_i')_{(0)},(\bar{a}_j')_{(0)}\rangle 
-\langle \bar{a}_i',\bar{a}_j'\rangle }\prod_{1 \le i\le k}\left(1+\frac{w_i}{z}\right)
^{-\langle (\bar{a}_i')_{(0)},(\bar{a}_i')_{(0)}\rangle /2
+\langle \bar{a}_i',\bar{a}_i'\rangle /2}.\nonumber
\end{eqnarray}
{}From  (\ref{main1}),
\begin{equation}\label{5.24}
A=\mbox{$\circ\atop\circ$}Y_*(a_1,w_1)\cdot\cdot\cdot 
Y_*(a_k,w_k)\mbox{$\circ\atop\circ$}\cdot \iota(1).
\end{equation}
The action of $e^{\Delta_z}$ on $A$ is given by:
\begin{eqnarray}
& &e^{\Delta_z}A=\prod_{i,j=1}^k \prod_{r=0}^{p-1} {\rm exp} 
\left(\sum_{s,t\ge 0, (s,t)\ne 0} c_{str}\langle \nu^r
\bar{a}'_i,\bar{a}'_j\rangle 
\left(\frac{w_i}{z}\right)^s\left(\frac{w_j}{z}\right)^t\right)A\nonumber\\
& &\ \ \ \ \ \ \ \  =\prod_{i,j=1}^k \prod_{r=0}^{p-1} {\rm exp}
\left(\sum_{s,t\ge 0, (s,t)\ne 0} c_{str}\left(\frac{w_i}{z}\right)^s\left(\frac{w_j}{z}\right)^t
\right)^{\langle \nu^r \bar{a}'_i,\bar{a}'_j\rangle }A\nonumber\\
& &\ \ \ \ \ \ \ \  =\prod_{i,j=1}^k \prod_{r=1}^{p-1} 
\left(\frac{(1+\frac{w_i}{z})^{1/p}-\omega^{-r}_p(1+\frac{w_j}{z})^{1/p}}
{1-\omega^{-r}_p}\right)^{-\langle \bar{a}'_i,\bar{a}'_j\rangle /2+\langle \nu^r\bar{a}'_i,\bar{a}'_j\rangle /2}A\label{5.17+}\\
& &\ \ \ \ \ \ \ \  =\prod_{1\le i<j\le k} \prod_{r=1}^{p-1} 
\left(\frac{(1+\frac{w_i}{z})^{1/p}-\omega^{-r}_p(1+\frac{w_j}{z})^{1/p}}
{1-\omega^{-r}_p}\right)^{-\langle
\bar{a}'_i,\bar{a}'_j\rangle+\langle
\nu^r\bar{a}'_i,\bar{a}'_j\rangle}\cdot\nonumber\\
& &\ \ \ \ \ \ \ \ \ \ \ \cdot \prod_{1\le i\le k}\left(1+\frac{w_i}{z}\right)^
{\langle (\bar{a}_i')_{(0)},(\bar{a}_i')_{(0)}\rangle /2
-\langle \bar{a}_i',\bar{a}_i'\rangle /2}A,\nonumber
\end{eqnarray}
where these expansions are to be understood as the formal power series
in the $w$'s that they came from. 
{}From (\ref{5.23})-(\ref{5.17+}) we find that
\begin{equation}\label{5.26}
Y^{\nu}_*(A,z)=\mbox{$\circ\atop\circ$}
Y^{\nu}_*(a_1,z+w_1)\cdot\cdot\cdot Y^{\nu}_*(a_k,z+w_k)
\mbox{$\circ\atop\circ$}f_{\bar a}(z,w_1,...,w_k)c_{\bar a_1,...\bar
a_k}\
\end{equation}
where
\begin{equation}\label{5.18+}
c_{\bar a_1,...,\bar a_k}=\frac{\sigma(\bar a)}{\sigma(\bar a_1)
\cdots\sigma(\bar a_k)}\prod_{1\le i<j\le k}\frac{\epsilon_2(\bar a_i,\bar a_j)}{\epsilon_1(\bar a_i,\bar a_j)}.
\end{equation}
(Lemma \ref{ltau} is used here.) Similarly,
\begin{equation}\label{5.27}
Y^{\nu}_*(B,z)=\mbox{$\circ\atop\circ$}
Y^{\nu}_*(b_1,z+x_1)\cdot\cdot\cdot Y^{\nu}_*(b_l,z+x_l)
\mbox{$\circ\atop\circ$}f_{\bar b}(z,x_1,...,x_l)c_{\bar b_1,...,\bar b_l}
\end{equation}
where $c_{\bar b_1,...,\bar b_l}$ is defined as in (\ref{5.18+}). 
Thus
\begin{eqnarray}
& & \mbox{$\circ\atop\circ$}
Y^{\nu}_*(A,z_1)Y^{\nu}_*(B,z_2)
\mbox{$\circ\atop\circ$}=f_{\bar a}(z_1,w_1,...,w_k)
f_{\bar b}(z_2,x_1,...,x_l)c_{\bar a_1,...,\bar a_k}
c_{\bar b_1,...,\bar b_l}\cdot\nonumber\\
& &\ \ \ \ \ \ \ \cdot\mbox{$\circ\atop\circ$}
Y^{\nu}_*(a_1,z_1+w_1)\cdot\cdot\cdot Y^{\nu}_*(a_k,z_1+w_k)Y^{\nu}_*(b_1,z_2+x_1)\cdot\cdot\cdot Y^{\nu}_*(b_l,z_2+x_l)
\mbox{$\circ\atop\circ$}\nonumber\\
& &\ \ \ \ \ \ \ \ \ \ \ \ \ \in ({\rm End }\,V^T_L)\{z_1,z_2\}[[w_1,\cdot\cdot\cdot,w_k,x_1,
\cdot\cdot\cdot,x_l]],\label{5.28}
\end{eqnarray}
and $\displaystyle{\lim_{z_1^{1/p}\to z_2^{1/p}}}\mbox{$\circ\atop\circ$}Y^{\nu}_*(A,z_1)Y^{\nu}_*(B,z_2)
\mbox{$\circ\atop\circ$}$ exists. That is, if 
$A(z_1,z_2)$ is the coefficient of any fixed monomial in $w_1,...,x_l$ in
$\mbox{$\circ\atop\circ$}Y^{\nu}_*(A,z_1)Y^{\nu}_*(B,z_2)\mbox{$\circ\atop\circ$},$ then  $\displaystyle{\lim_{z_1^{1/p}\to z_2^{1/p}}}A(z_1,z_2)$ exists.
{}From (\ref{9.20}) we see that
\begin{equation}\label{5.29}
Y^{\nu}_*(A,z_1)Y^{\nu}_*(B,z_2)=
\mbox{$\circ\atop\circ$}Y^{\nu}_*(A,z_1)Y^{\nu}_*(B,z_2)\mbox
{$\circ\atop\circ$}f_{(\bar a,\bar b)}(z_1,w_1,...,w_k;z_2,x_1,...,x_l).
\end{equation}

Fix a monomial
\begin{equation}\label{5.30}
P=\prod_{1\le i\le k, 1\le j\le l}w^{r_i}_ix^{s_j}_j (r_i, s_j\ge 1)
\end{equation}
in the $w_i$ and $x_j$. We may and do choose $N\ge 0$ so large that the
coefficients of $P$ and of each monomial of lower total degree than
$P$ in
$$F_N=(z_1-z_2)^{N}
f_{(\bar a,\bar b)}(z_1,w_1,...,w_k;z_2,x_1,...,x_l)F_{(\bar a,\bar b)}(z_1,z_2)$$
are polynomials in $z_1^{-1},$ $z_2^{-1},$ $z_1^{1/p}$ and $z_2^{1/p}.$
In fact, the coefficient of $P$ in $F_N$  is 
\begin{eqnarray*}
& &\prod_{1\leq i\leq k}\prod_{1\leq j\leq l}\left(\frac{\partial}{\partial w_i}\right)^{r_i}\left(\frac{\partial}{\partial x_j}\right)^{s_j}F_N|_{w_i=x_j=0}\\
& &\ \ =\sum C_{mnI}\prod_{r=0}^{p-1}(z_1^{1/p}-\omega_p^{-r}z_2^{1/p})^{i_r}z_1^{1/p-m}z_2^{1/p-n},
\end{eqnarray*}
where the sum is over $m,$ $n,$ $I=(i_0,...,i_{p-1})$ with $0\leq i_r,
m, n\leq N,$ and the 
$C_{mnI}$ are constants. Then the coefficients of $P$ and of each monomial 
of lower total degree than $P$ in 
$$(-z_2+z_1)^{N}
f_{(\bar b,\bar a)}(z_2,x_1,...,x_l;z_1,w_1,...,w_k)F_{(\bar b,\bar a)}(z_2,z_1)$$
are also  polynomials in $z_i^{-1}$ and $z_i^{1/p}$ for $i=1,2,$ 
and agree with the corresponding coefficients in $F_N.$ 
Let $Y_P(z_1^{1/p},z_2^{1/p})$ denote the coefficient of $P$ in
\begin{equation}\label{5.31}
Y^{\nu}_*(A,z_1)Y^{\nu}_*(B,z_2)(z_1-z_2)^NF_{(\bar a,\bar b)}(z_1,z_2)
=\mbox{$\circ\atop\circ$}Y^{\nu}_*(A,z_1)Y^{\nu}_*(B,z_2)\mbox
{$\circ\atop\circ$}F_N.
\end{equation}
By ({\ref{5.28}) and (\ref{5.29}) we see that
\begin{equation}\label{5.32}
\begin{array}{c}
Y_P(z_1^{1/p},z_2^{1/p})z_1^{-\bar{a}'_{(0)}-\langle \bar{a}'_{(0)},\bar{a}'_{(0)}\rangle 
/2+\langle \bar{a}',\bar{a}'\rangle /2}z_2^{-\bar{b}'_{(0)}-\langle \bar{b}'_{(0)},\bar{b}'_{(0)}\rangle 
/2+\langle \bar{b}',\bar{b}'\rangle /2}\\
\in ({\rm End}\,V^T_L)[[z_1^{1/p},z_1^{-1/p},z_2^{1/p},z_2^{-1/p}]],
\end{array}
\end{equation}
and ${\displaystyle \lim_{z_1^{1/p}\to z_2^{1/p}}Y_P(z_1^{1/p},z_2^{1/p})}$
exists. The coefficient of $P$ in
$$Y^{\nu}_*(A,z_1)Y^{\nu}_*(B,z_2)(z_1-z_2)^nF_{(\bar a,\bar b)}(z_1,z_2)$$
is
$$Y_P(z_1^{1/p},z_2^{1/p})(z_1-z_2)^{-(N-n)}=Y_P(z_1^{1/p},z_2^{1/p})z_1^{-(N-n)}\left(1-\frac{z_2}{z_1}\right)
^{-(N-n)}.$$

Similarly, reversing the roles of $A$ and $B$ and of $z_1$ and $z_2$
and noting that
\begin{equation}\label{5.33}
\mbox{$\circ\atop\circ$}Y^{\nu}_*(B,z_2)Y^{\nu}_*(A,z_1)
\mbox{$\circ\atop\circ$}=(c_{\nu}(\bar{a},\bar{b}))^{-1}\mbox{$\circ\atop\circ$}
Y^{\nu}_*(A,z_1)Y^{\nu}_*(B,z_2)\mbox{$\circ\atop\circ$},
\end{equation}
we find that
$$Y^{\nu}_*(B,z_2)Y^{\nu}_*(A,z_1)=c_{\nu}(\bar{a},\bar{b})^{-1}
\mbox{$\circ\atop\circ$}Y^{\nu}_*(A,z_1)Y^{\nu}_*(B,z_2)
\mbox{$\circ\atop\circ$}f_{(\bar b,\bar a)}(z_2,x_1,...,x_l;z_1,w_1,...,w_k)$$
and that the coefficient of $P$ in
$$c_{\nu}(\bar{a},\bar{b})Y^{\nu}_*(B,z_2)Y^{\nu}_*(A,z_1)(-z_2+z_1)^NF_{(\bar b,\bar a)}(z_2,z_1)$$
is also $Y_P(z_1^{1/p},z_2^{1/p}).$ Thus the coefficient of $P$ in
$$c_{\nu}(\bar{a},\bar{b})Y^{\nu}_*(B,z_2)Y^{\nu}_*(A,z_1)(-z_2+z_1)^n
F_{(\bar b,\bar a)}(z_2,z_1)$$
is
$$Y_P(z_1^{1/p},z_2^{1/p})(-z_2+z_1)^{-(N-n)}=Y_P(z_1^{1/p},z_2^{1/p})(-z_2)^{-(N-n)}\left(1-\frac{z_1}
{z_2}\right)^{-(N-n)}.$$
Applying Proposition \ref{10.1} with $A(z_1^{1/p},z_2^{1/p})$ equal to 
$Y_P(z_1^{1/p},z_2^{1/p}),$ we see that the coefficient of $P$ in
\begin{equation}\label{5.34}
\begin{array}{l}
F_{(\bar a,\bar b)}(z_1,z_2)
(z_1-z_2)^nY^{\nu}_*(A,z_1)Y^{\nu}_*(B,z_2)\\
\ \ -c_{\nu}(\bar{a},\bar{b})F_{(\bar b,\bar a)}(z_2,z_1)
(-z_2+z_1)^nY^{\nu}_*(B,z_2)Y^{\nu}_*(A,z_1)
\end{array}
\end{equation}
is
\begin{equation}\label{5.35} 
\begin{array}{l}
\displaystyle{
\mbox{Res}_{z_0}z_2^{-1}z_0^{n-N}
Y_P(z_1^{1/p},z_2^{1/p})e^{-z_0\frac{\partial}{\partial
z_1}}\delta\left(\frac{z_1}{z_2}\right)} \\
\displaystyle{
\ \ \ \ =\mbox{Res}_{z_0}\frac{1}{p}z_2^{-1}z_0^{n-N}\sum_{r=0}^{p-1}
Y_P(z_1^{1/p},z_2^{1/p})e^{-z_0\frac{\partial}{\partial
 z_1}}\delta\left(\omega^r_p\left(\frac{z_1}{z_2}\right)^{1/p}\right)}
 \\
\displaystyle{
\ \ \ \
=\mbox{Res}_{z_0}\frac{1}{p}z_2^{-1}z_0^{n-N}\sum_{r=0}^{p-1}
\left(\left(\lim_{z_1^{1/p}\to
\omega^{-r}_pz_2^{1/p}}Y_P((z_1+z_0)^{1/p},z_2^{1/p})\right)\cdot\right.} \\
\displaystyle{
\ \ \ \ \  \ \ \  \cdot \omega_p^{r\sum_s\nu^s\bar{a}'+r\langle \sum_s\nu^s\bar{a}',\bar{a}'\rangle /2-pr\langle \bar{a}',\bar{a}'\rangle /2}\cdot} \\
\displaystyle{
\ \ \ \ \ \ \ \ \left.\cdot e^{-z_0\frac{\partial}{\partial z_1}}\left(\frac{z_1}{z_2}\right)^
{\bar{a}'_{(0)}+\langle \bar{a}'_{(0)},\bar{a}'_{(0)}\rangle /2
-\langle \bar{a}',\bar{a}'\rangle /2}\delta\left(\omega^r_p\left(\frac{z_1}{z_2}\right)^{1/p}\right)\right),}
\end{array}
\end{equation}
where $s$ ranges over ${\Bbb Z}/p{\Bbb Z}.$ Now we
compute
$\displaystyle{\lim_{z_1^{1/p}\to\omega^{-r}_pz_2^{1/p}}}Y_P((z_1+z_0)^{1/p},z_2^{1/p}).$
Note that $Y_P((z_1+z_0)^{1/p},z_2^{1/p})$ is the coefficient of $P$ in
$$\mbox{$\circ\atop\circ$}Y^{\nu}_*(A,z_1+z_0)Y^{\nu}_*(B,z_2)\mbox
{$\circ\atop\circ$}(z_1+z_0-z_2)^NF_{(\bar a,\bar b)}(z_1+z_0,z_2)
f_{(\bar a,\bar b)}(z_1+z_0,w_1,...,w_k;z_2,x_1,...,x_l).$$ From
(\ref{extral}) we have
$$\lim_{z_1^{1/p}\to\omega_p^{-r}z_2^{1/p}}\mbox{$\circ\atop\circ$}
Y^{\nu}_*(A,z_1+z_0)Y^{\nu}_*(B,z_2)\mbox{$\circ\atop\circ$}$$
$$=\mbox{$\circ\atop\circ$}Y^{\nu}_*(\hat{\nu}^rA,z_2+z_0)
Y^{\nu}_*(B,z_2)\mbox{$\circ\atop\circ$}c_{\nu}(\bar a-\nu^r\bar
a,\bar b)(\hat\nu^ra^{-1})a
\omega_p^{-r\sum\nu^s\bar{a}'-r\langle \sum\nu^s\bar{a}',\bar{a}'\rangle /2
+rp\langle \bar{a}',\bar{a}'\rangle /2}.$$
We also have
$$\lim_{z_1^{1/p}\to \omega_p^{-r}z_2^{1/p}}F_{(\bar a,\bar b)}(z_1+z_0,z_2)
f_{(\bar a,\bar b)}(z_1+z_0,w_1,...,w_k;z_2,x_1,...,x_l)$$
$$=F_{(\nu^r\bar a,\bar b)}(z_2+z_0,z_2)f_{(\nu^r\bar a,\bar b)}
(z_2+z_0,w_1,...,w_k;z_2,x_1,...,x_l).$$
Thus $\displaystyle{\lim_{z_1^{1/p}\to\omega_p^{-r}z_2^{1/p}}}Y_P((z_1+z_0)^{1/p},z_2^{1/p})$ is the
coefficient of $P$ in
\begin{equation}\label{5.36}
\begin{array}{c}
\mbox{$\circ\atop\circ$}Y^{\nu}_*(\hat{\nu}^rA,z_2+z_0)
Y^{\nu}_*(B,z_2)\mbox{$\circ\atop\circ$}c_{\nu}(\bar
a-\nu^r\bar a,\bar b)(\hat\nu^ra^{-1})az_0^NF_{(\nu^r\bar a,\bar b)}(z_2+z_0,z_2)\cdot\\
\cdot f_{(\nu^r\bar a,\bar b)}
(z_2+z_0,w_1,...,w_k;z_2,x_1,...,x_l)
\omega_p^{-r\sum\nu^s\bar{a}'-r\langle \sum\nu^s\bar{a}',\bar{a}'\rangle /2
+rp\langle \bar{a}',\bar{a}'\rangle /2}.
\end{array}
\end{equation}
We conclude that the coefficient of $P$ in (\ref{5.34}) is the coefficient
of $P$ in
$$\mbox{Res}_{z_0}\frac{1}{p}z_2^{-1}z_0^{n}\sum_{r=0}^{p-1}\left(
\mbox{$\circ\atop\circ$}Y^{\nu}_*(\hat{\nu}^rA,z_2+z_0)
Y^{\nu}_*(B,z_2)\mbox{$\circ\atop\circ$}
\delta\left(\omega^r_p\left(\frac{z_1-z_0}{z_2}\right)^{1/p}\right)\cdot
\right.$$
\begin{equation}\label{5.37}
\cdot F_{(\nu^r\bar a,\bar b)}(z_2+z_0,z_2)f_{(\nu^r\bar a,\bar b)}
(z_2+z_0,w_1,...,w_k;z_2,x_1,...,x_l)
\cdot
\end{equation}
$$
\left.\cdot c_{\nu}(\bar a-\nu^r\bar a,\bar b)(\hat\nu^{r}a^{-1})a\left(\frac{z_1-z_0}{z_2}\right)^
{\bar{a}'_{(0)}+\langle \bar{a}'_{(0)},\bar{a}'_{(0)}\rangle /2
-\langle \bar{a}',\bar{a}'\rangle /2}
\right).$$
Therefore the operators in ({\ref{5.34}) and in (\ref{5.37}) are the same.
This last assertion is independent of $P$ and $N$.

On the other hand, just as in the proof of Theorem 5.1 of [DL2], we obtain
$$Y_*(\hat{\nu}^rA,z_0)B=
\mbox{$\circ\atop\circ$}
Y_*(\hat{\nu}^ra_1,z_0+w_1)\cdot\cdot\cdot Y_*(\hat\nu^ra_k,z_0+w_k)
Y_*(b_1,x_1)\cdot\cdot\cdot Y_*(b_l,x_l)
\mbox{$\circ\atop\circ$}\iota(1)\cdot$$
$$\cdot\prod_{1\le i\le k,\, 1\le j\le l}
(z_0+w_i-x_j)^{\langle \hat{\nu}^r\bar{a}'_i,\bar{b}'_j\rangle }.$$
By (\ref{5.24}) and (\ref{5.26}) we have
\begin{equation}\label{5.38}
\begin{array}{l}
Y^{\nu}_*(Y_*(\hat{\nu}^rA,z_0)B,z_2)=\mbox{$\circ\atop\circ$}
Y^{\nu}_*(\hat\nu^ra_1,z_2+z_0+w_1)\cdot\cdot\cdot
Y^{\nu}_*(\hat\nu^ra_k,z_2+z_0+w_k)\cdot\\
\hspace*{2 cm}\cdot Y^{\nu}_*(b_1,z_2+x_1)\cdot\cdot\cdot 
Y^{\nu}_*(b_l,z_2+x_l)\mbox{$\circ\atop\circ$}
f_{\bar a}(z_2+z_0,w_1,...,w_k)\cdot\\
\hspace*{2 cm}\cdot f_{\bar b}(z_2,x_1,...,x_l)f_{(\nu^r\bar a,\bar b)}
(z_2+z_0,w_1,...,w_k;z_2,x_1,...,x_l)\cdot\\
\hspace{2 cm} 
\displaystyle{
\cdot\prod_{1\le i\le k,\, 1\le j\le l}} 
\, {\displaystyle{\prod_{s=0}^{p-1}}\left((z_2+z_0+w_i)^{1/p}-
\omega_p^s(z_2+x_j)^{1/p}\right)^{-\langle \nu^r\bar{a}'_i,
\bar{b}'_j\rangle}}\cdot\\
\hspace*{2 cm}\cdot\displaystyle{\prod_{1\le i\le k,\, 1\le j\le l}(z_0+w_i-x_j)^{\langle \hat{\nu}^r\bar{a}'_i,\bar{b}'_j\rangle }}c_{\nu^r\bar a_1,...,\nu^r\bar a_k}c_{\bar b_1,...,\bar b_l}
(\tau(\nu^r\bar a,\bar b))^{-1}
\end{array}
\end{equation}
where we have used the relation
$$c_{\nu^r\bar a_1,...,\nu^r\bar a_k,\bar b_1,...,\bar b_l}=
c_{\nu^r\bar a_1,...,\nu^r\bar a_k}c_{\bar b_1,...,\bar b_l}
(\tau(\nu^r\bar a,\bar b))^{-1}$$
(see (\ref{e5.8}) and (\ref{5.18+})). We next prove that
$$\displaystyle{
\prod_{1\le i\le k,\, 1\le j\le l}}\, 
{\displaystyle{\prod_{s=0}^{p-1}}\left((z_2+z_0+w_i)^{1/p}-
\omega_p^s(z_2+x_j)^{1/p}\right)^{-\langle \nu^r\bar{a}'_i,\bar{b}'_j\rangle }}\displaystyle{\prod_{1\le i\le k,\, 1\le j\le l}(z_0+w_i-x_j)^{\langle \hat{\nu}^r\bar{a}'_i,\bar{b}'_j\rangle }}$$
$$=z_0^{\langle\nu^r\bar a',\bar b'\rangle}
\prod_{s=0}^{p-1}\left((z_2+z_0)^{1/p}-
\omega_p^sz_2^{1/p}\right)^{-\langle\nu^r\bar a',\bar
b'\rangle}.$$
For this purpose we introduce
\begin{equation}\label{e5.33}
g(z_2,z_0,w,x)=\prod_{s=0}^{p-1}\left((z_2+z_0+w)^{1/p}-
\omega_p^s(z_2+x)^{1/p}\right)^c(z_0+w-x)^{-c}
\end{equation}
where $c\in {\Bbb C},$ $w,x$ are formal variables, and the formal 
series is interpreted as above. It suffices to prove
that $g(z_2,z_0,w,x)=g(z_2,z_0,0,0),$ or that 
$\frac{\partial g}{\partial w}=\frac{\partial g}{\partial x}=0.$ Now we calculate $\frac{\partial g}{\partial w}:$
$$\frac{\partial g}{\partial w}=\sum_{s=0}^{p-1}\frac{c}{p}
(z_2+z_0+w)^{1/p-1}\left((z_2+z_0+w)^{1/p}-
\omega_p^s(z_2+x)^{1/p}\right)^{-1}g-c(z_0+w-x)^{-1}g,$$
which is zero because
$$\sum_{s=0}^{p-1}\frac{1}{p}
(z_2+z_0+w)^{1/p-1}\left((z_2+z_0+w)^{1/p}-
\omega_p^s(z_2+x)^{1/p}\right)^{-1}=(z_0+w-x)^{-1}.$$
Similarly, $\frac{\partial g}{\partial x}=0.$ 

Finally, we see from (5.28) that the operator 
$Y_*^{\nu}(Y_*(\hat\nu^rA,z_0)B,z_2)$ is 
equal to 
\begin{equation}\label{extral1}
\begin{array}{c}
\mbox{$\circ\atop\circ$}
Y^{\nu}_*(\hat{\nu}^rA,z_2+z_0)Y^{\nu}_*(B,z_2)
\mbox{$\circ\atop\circ$}
f_{({\nu}^r\bar a,\bar b)}(z_2+z_0,w_1,...,w_k;z_2,x_1,...,x_l)\cdot\\
\hspace*{2 cm}
\displaystyle{
\cdot (\tau(\nu^r\bar a,\bar b))^{-1}z_0^{\langle\nu^r\bar a',\bar b'\rangle}
\prod_{s=0}^{p-1}\left((z_2+z_0)^{1/p}-
\omega_p^sz_2^{1/p}\right)^{-\langle\nu^r\bar a',\bar
b'\rangle}}.
\end{array}
\end{equation}
Combining (\ref{5.34}), (\ref{5.35}), (\ref{5.36}), (\ref{5.37}) and
(\ref{extral1}) completes the proof.
\ \ \ \ \ \ \ $\Box$   

\begin{rem}\label{rem.5.1}\hspace*{-0.2 cm}{\bf : }{\rm In the case in which $L$
is a direct sum of several copies of the root lattice of the Lie 
algebra $sl(2,{\Bbb C})$ and $\nu=-1,$ Theorem \ref{t10.2} was obtained
in [Hu], and both the twisted $Z$-algebra relations [LW1]-[LW5]
and the twisted parafermion algebra relations [ZF2] associated with
the twisted vertex operator constructions of $A_1^{(1)}$ were recovered as a
consequence of this Jacobi identity  and a 
multi-operator extension of it.}
\end{rem}

The Jacobi  identity  for relative untwisted vertex operators 
([DL1], [DL2], [FFR]) is a
special case of Theorem \ref{t10.2}: 
\begin{coro}\label{c.5.1}
Let $u,v\in V_L.$ Then 
$$z^{-1}_0\left(\frac{z_1-z_2}{z_0}\right)^{-\langle\bar a',\bar b'\rangle}
\delta\left(\frac{z_1-z_2}{z_0}\right)
Y_*(u,z_1)Y_*(v,z_2)$$
\begin{equation}\label{ujac}
-c(\bar{a},\bar{b})
z^{-1}_0\left(\frac{z_2-z_1}{z_0}\right)^{-\langle\bar a',\bar b'\rangle}
\delta\left(\frac{z_2-z_1}{-z_0}\right)
Y_*(v,z_2)Y_*(u,z_1)
\end{equation}
$$=z_2^{-1}\delta\left(\frac{z_1-z_0}{z_2}\right)Y_*(Y_*(u,z_0)v,z_2)
\left(\frac{z_1-z_0}{z_2}\right)^{\bar a'}.$$
\end{coro}

\noindent{\bf Proof\ } If $\nu=1,\ p=1,\ c^{\nu}_0(\cdot,\cdot)=
c_0(\cdot,\cdot)$ and $T={\Bbb C}\{L\},$ then $V_L^T=V_L$ and
$Y_*^{\nu}(u,z)=Y_*(u,z).$ Moreover, $F_{(\alpha ,\beta)}=
(z_1-z_2)^{-\langle\alpha',\beta'\rangle},$ $G_{(\alpha,\beta)}=1$ and
$\tau(\alpha,\beta)=1$ for $\alpha,\beta\in L$ in this case.
Now the Jacobi identity (\ref{ujac}) follows 
immediately from (\ref{tjac}).
\ \ \ \ \ $\Box$  

We also include the Jacobi identity for ordinary (unrelativized) twisted vertex
operators (and therefore the commutator formula for these
operators ([FLM2], [L1])) as follows:
\begin{coro}\label{co1} In the settings of
Remarks 2.2, \ref{erem4.4} and \ref{erem4.2},
for $u, v\in V_L,$ we have 
$Y^{\nu}_*(u,z)=Y_{\nu}(u,z)$ and 
\begin{equation}\label{tjac2}
\begin{array}{c}
\displaystyle{
z^{-1}_0
\delta\left(\frac{z_1-z_2}{z_0}\right)Y_{\nu}(u,z_1)Y_{\nu}(v,z_2)-
z^{-1}_0
\delta\left(\frac{z_2-z_1}{-z_0}\right)Y_{\nu}(v,z_2)Y_{\nu}(u,z_1)}\\
\displaystyle{
=\frac{1}{p}z^{-1}_2\displaystyle{\sum^{p-1}_{r=0}}\delta\left(\omega^r_p\left(\frac{z_1-z_0}{z_2}\right)^{1/p}\right)Y_{\nu}(Y(\hat\nu^r
u,z_0)v,z_2).}
\end{array}
\end{equation}
\end{coro}

\noindent {\bf Proof } From (\ref{add1}) and (\ref{right}),
$$c_{\nu}(\bar a-\nu^r\bar a,\bar b)=\prod_{s=0}^{p-1}(-\omega_p^s)^{\langle \nu^s(\bar
a-\nu^r\bar a),\bar b\rangle}=\omega_p^{r\langle\sum\nu^s\bar a,\bar b\rangle}$$
$$(\hat\nu^ra^{-1})a=\omega_p^{r\sum\nu^s\bar a+r\langle\sum\nu^s\bar
a,\bar a\rangle/2}.$$
Since 
$$\langle \sum_s\nu^s\bar{a},\bar{a}\rangle  \in 2{\Bbb Z} \ \ {\rm
for}\ \ \alpha\in L$$
by (\ref{L2.1}) we see that 
$$\langle \bar{a}_{(0)},\bar{a}_{(0)}\rangle \in \frac{2}{p}{\Bbb Z}$$
and that
$$\left(\frac{z_1-z_0}{z_2}\right)^{\bar{a}_{(0)}}\in
({\rm End}\,V^T_L)[[z_1^{1/p},z_1^{-1/p},z_2^{1/p},z_2^{-1/p},z_0]].$$
Using the basic property of the $\delta$-function we have:
$$\left(\frac{z_1-z_0}{z_2}\right)^{\bar{a}_{(0)}+\langle \bar{a}_{(0)},\bar{a}_{(0)}\rangle /2
-\langle \bar{a},\bar{a}\rangle /2}\delta\left(\omega^r_p\left(\frac{z_1-z_0}{z_2}\right)^{1/p}\right)$$
$$=
\omega_p^{-r\sum\nu^s\bar{a}-r\langle \sum\nu^s\bar{a},\bar{a}\rangle /2}
\delta\left(\omega^r_p\left(\frac{z_1-z_0}{z_2}\right)^{1/p}\right).$$
Since $L$ is even we see that
$$G_{(\nu^r\bar a,\bar b)}(z_2+z_0,z_2)z_0^{-\langle\nu^r\bar a,\bar b\rangle}=F_{(\bar a,\bar b)}(z_2+z_0,z_2).$$
Then the Jacobi identity (\ref{tjac}) reduces to:
$$z_0^{-1}F_{(\bar a,\bar
b)}(z_1,z_2)\delta\left(\frac{z_1-z_2}{z_0}\right)
Y_{\nu}(u,z_1)Y_{\nu}(v,z_2)$$
\begin{equation}\label{tjac1}
-c_{\nu}(\bar{a},\bar{b})z_0^{-1}
F_{(\bar b,\bar a)}(z_2,z_1)\delta\left(\frac{z_2-z_1}{-z_0}\right)Y_{\nu}(v,z_2)Y_{\nu}(u,z_1)
\end{equation}
$$=\frac{1}{p}
z_2^{-1}\sum_{r=0}^{p-1}\left(F_{(\nu^r\bar a,\bar b)}(z_2+z_0,z_2)
Y_{\nu}(Y(\hat{\nu}^ru,z_0)v,z_2)\omega_p^{r\langle \sum\nu^s\bar{a},\bar{b}\rangle }
\delta\left(\omega^r_p\left(\frac{z_1-z_0}{z_2}\right)^{1/p}\right)\right)$$
(recall that $\tau(\alpha,\beta)=1$ for $\alpha,\beta\in L$ in this
case). Noting that 
$$\prod_{s=0}^{p-1}(z_2^{1/p}-\omega^{-s}_pz_1^{1/p})^{-\langle \nu^s\bar{b},
\bar{a}\rangle }=(c_{\nu}(\bar{a},\bar{b}))^{-1}\prod_{s=0}^{p-1}(-\omega^{-s}_pz_2^{1/p}
+z_1^{1/p})^{-\langle \nu^s\bar{a},\bar{b}\rangle }$$
and that at least one of the two binomials $(z_1^{1/p}-\omega_p^sz_2^{1/p})^
{\langle \nu^s\bar{a},\bar{b}\rangle }$ and $(-\omega_p^sz_2^{1/p}+z_1^{1/p})^
{-\langle \nu^s\bar{a},\bar{b}\rangle }$ is a polynomial in $z_1^{1/p}$
and $z_2^{1/p},$ and multiplying (\ref{tjac1}) by
$\displaystyle{\prod_{s=0}^{p-1}}(z_1^{1/p}-\omega^{-s}_pz_2^{1/p})^{\langle \nu^s\bar{a},\bar{b}\rangle },$ 
we have:
$$z_0^{-1}\delta\left(\frac{z_1-z_2}{z_0}\right)Y_{\nu}(u,z_1)Y_{\nu}(v,z_2)-
z_0^{-1}\delta\left(\frac{z_2-z_1}{-z_0}\right)Y_{\nu}(v,z_2)Y_{\nu}(u,z_1)$$
$$=\frac{1}{p}
z_2^{-1}\sum_{r=0}^{p-1}\left(F_{(\nu^r\bar a,\bar b)}(z_2+z_0,z_2)
\prod_s(z_1^{1/p}-\omega^{-s}_pz_2^{1/p})^
{\langle \nu^s\bar{a},\bar{b}\rangle }\right.\cdot$$
$$\left.\cdot Y_{\nu}(Y(\hat{\nu}^ru,z_0)v,z_2)\omega_p^{r\langle \sum\nu^s\bar{a},\bar{b}\rangle }
\delta\left(\omega^r_p\left(\frac{z_1-z_0}{z_2}\right)^{1/p}\right)\right).$$
Now (\ref{tjac2}) follows from: 
$$\prod_{s=0}^{p-1}(z_1^{1/p}-\omega^{-s}_pz_2^{1/p})^
{\langle \nu^s\bar{a},\bar{b}\rangle }
\delta\left(\omega^r_p\left(\frac{z_1-z_0}{z_2}\right)^{1/p}\right)
\hspace*{1 cm}$$
$$=(F_{(\nu^r\bar a,\bar b)}(z_2+z_0,z_2))^{-1}\omega^{-r\langle \sum\nu^s\bar{a},\bar{b}\rangle }
\delta\left(\omega^r_p\left(\frac{z_1-z_0}{z_2}\right)^{1/p}\right).\ \ \ \ \ \ \Box$$
 
The following corollary is used to construct the ``moonshine
modules'' based on order $p$ isometries of the Leech lattice for the odd
primes 3, 5, 7 and 13 [DM1]:
\begin{coro}\label{co2}
Let $p$ be an odd prime and let $\nu$ be an
isometry of $L$ having no nonzero fixed points. 
Let $L_0$ be an even sublattice of $L$ such that
(\ref{L2.2}), (\ref{add1}) and (\ref{LL.2}) in Remark 2.2 hold for the restrictions of
$c_0(\cdot,\cdot),$ $c_0^{\nu}(\cdot,\cdot)$ and
$\epsilon_0(\cdot,\cdot)$ to $L_0.$  Take ${\bf h}_*=0.$
Then  for $u, v\in V_{L_0},$ we have $Y^{\nu}_*(u,z)=Y_{\nu}(u,z)$ and 
\begin{equation}\label{cjac2}
\begin{array}{c}
\displaystyle{
z^{-1}_0
\delta\left(\frac{z_1-z_2}{z_0}\right)Y_{\nu}(u,z_1)Y_{\nu}(v,z_2)-
z^{-1}_0
\delta\left(\frac{z_2-z_1}{-z_0}\right)Y_{\nu}(v,z_2)Y_{\nu}(u,z_1)}\\
\displaystyle{
=\frac{1}{p}z^{-1}_2\displaystyle{\sum^{p-1}_{r=0}}\delta\left(\omega^r_p\left(\frac{z_1-z_0}{z_2}\right)^{1/p}\right)Y_{\nu}(Y(\hat\nu^ru,z_0)v,z_2)
(\hat\nu^ra^{-1})a}
\end{array}
\end{equation}
as operators on $V_L^T.$
\end{coro}

\noindent {\bf Proof }
Note that  ${\bf h}_{(0)}=0$ and $\alpha_{(0)}=0$ for $\alpha\in {\bf h}$
because $\nu$ is fixed-point free. Now the proof is almost the same as that of
Corollary \ref{co1}, except that we do not 
replace $(\hat\nu^ra^{-1})a$ by $\omega_p^{r\sum\nu^s\bar a+r\langle\sum\nu^s\bar
a,\bar a\rangle/2
}.$\ \ \ \ \ $\Box$

\section{The Virasoro algebras}

\setcounter{equation}{0}

In this section we study the representations of the Virasoro algebra, 
with basis $\{L_n|n\in{\Bbb Z}\}\cup\{c\}$ and with 
the usual commutation relations 
\begin{equation}
[L_m,L_n]=(m-n)L_{m+n}+\frac{1}{12}(m^3-m)\delta_{m+n ,0}c
\ \mbox{for}\ \ m,n \in{\Bbb Z},
\end{equation}
on both $V_L$ and $V_L^T,$  in terms of the relative (untwisted and twisted)
vertex operators associated
with a canonical quadratic element of weight 2 in $V_L.$ 
The natural operator $\Delta_z$ (\ref{9.16}) incorporated into the
definition of the relative twisted vertex operators plays a
fundamental role. 
We use the action
of $L_0$ to reinterpret  the weight gradations of $V_L$ and $V_L^T.$
The reader can refer to [DL2] and [FLM3] for similar
discussions for relative untwisted vertex operators and background.

Recall that
$\{\beta_1,...,\beta_d\}$ is an orthonormal basis of
${\bf h}_*^{\perp}.$ 
Set
\begin{equation}\label{omega}
\omega=\frac{1}{2}\displaystyle{\sum^d_{i=1}}\beta_i(-1)^2.
\end{equation}
and set
\begin{equation}
L(z)=Y_*(\omega,z)=\displaystyle{\sum_{n\in{\Bbb Z}}}L(n)z^{-n-2},
\end{equation}
i.e,
\begin{equation} 
L(n)=\omega_{n+1}\ \ \mbox{for} \ \ n\in{\Bbb Z}
\end{equation}
(recall (\ref{3.41}), the definition of the components of a relative 
untwisted vertex operator). We also define operators $L^{\nu}(n)$ on $V^T_L$ as the
following generating function coefficients:
\begin{equation}
L^{\nu}(z)=Y^{\nu}_*(\omega,z)=\displaystyle{\sum_{n\in{\bf
Z}}}L^{\nu}(n)z^{-n-2},
\end{equation}
i.e,
\begin{equation} 
L^{\nu}(n)=\omega^{\nu}_{n+1}\ \ \mbox{for} \ \ n\in{\Bbb Z}
\end{equation}
(see (\ref{9.18})). 
 
We have the important relation between $L(-1)$ and
differentiation:
\begin{prop}\label{prop9}
For all $v\in V_L$,
\begin{equation}\label{7.7}
Y^{\nu}_*(L(-1)v,z)=\frac{d}{dz}Y^{\nu}_*(v,z).
\end{equation}
\end{prop}

\noindent {\bf Proof\ } First we observe that a special case of
(\ref{7.7}),
 for $v=\iota(a),$
follows from (\ref{9.21}) and the fact that
$$L(-1)\iota(a)=\bar{a}'(-1)\iota(a).$$
It is easy to
see from the definition (\ref{3.40}) that for all $v\in V_L,$
$$Y_*(L(-1)v,z)=\frac{d}{dz}Y_*(v,z).$$
Thus we have:
$$Y_*(e^{z_0L(-1)}v,z)=e^{z_0\frac{d}{dz}}Y_*(v,z)=Y_*(v,z+z_0).$$
We can now apply both sides to $\iota(1)$ and invoke (\ref{per1}) to obtain:
$$e^{z_0L(-1)}v=Y_*(v,z_0)\cdot\iota(1).$$
On the other hand, as operators on $V^T_L,$
$$Y_*^{\nu}(Y_*(v,z_0)\cdot\iota(1),z)=Y^{\nu}_*(v,z+z_0)$$
by the proof of Theorem \ref{t10.2}. Hence
$$e^{z_0\frac{d}{dz}}Y^{\nu}_*(v,z)=Y^{\nu}_*(v,z+z_0)=Y^{\nu}_*(e^{z_0L(-1)}v,z),$$
and ({\ref{7.7}) follows by extracting the coefficient of $z_0$.\ \ \ \ \ \ $\Box$

Combining the Jacobi identity with Proposition \ref{prop9} we have:
$$[L^{\nu}(z_1),Y^{\nu}_*(v,z_2)]=\mbox{Res}_{z_0}z_2^{-1}Y^{\nu}_*(L(z_0)v,z_2)
e^{-z_0\frac{\partial}{\partial z_1}}\delta(z_1/z_2)$$
\begin{equation}\label{7.8}
=z^{-1}_2\left(\frac{d}{dz}Y^{\nu}_*(v,z_2)\right)\delta(z_1/z_2)-z^{-1}_2
Y^{\nu}_*(L(0)v,z_2)\frac{\partial}{\partial z_1}\delta(z_1/z_2)
\end{equation}
$$+z^{-1}_2\mbox{Res}_{z_0}
\displaystyle{\sum_{n>0}}Y^{\nu}_*(L(n)v,z_2)z^{-n-2}_0e^{-z_0
\frac{\partial}{\partial z_1}}\delta(z_1/z_2)$$
for all $v\in V_L$. Equating the coefficient of $z^{-1}_1$ and changing
$z_2$ to $z$, we get
\begin{prop}
For all $v\in V_L$,
\begin{equation}
[L^{\nu}(-1),Y^{\nu}_*(v,z)]=\frac{d}{dz}Y^{\nu}_*(v,z)=Y^{\nu}_*(L(-1)v,z).\
\ \  \Box
\end{equation}
\end{prop}

We call a nonzero vector $v\in V_L$ (resp. $v\in V_L^T$)  a $weight$ $vector$ if $v$ satisfies the
following condition:
\begin{equation}\label{wt}
L(0)v=hv\ \ ({\rm resp.,\ } L^{\nu}(0)v=hv) \ \ \mbox{for some }\ \  
h\in{\Bbb C}
\end{equation}
and we call $h$ the weight of $v$. If $v$ further satisfies the condition
\begin{equation}
L(n)v=0\ \ ({\rm resp.,\ } L^{\nu}(n)v=0)\ \ \mbox{for}\ \ n>0, 
\end{equation}
we call $v$ a $lowest$ $weight$ $vector$.  

By (\ref{7.8}) we have:
\begin{prop}
If $v\in V_L$ is a lowest weight vector
with the weight $h$, then
\begin{equation}\label{7.11}
[L^{\nu}(z_1),Y^{\nu}_*(v,z_2)]=z_2^{-1}\left(\frac{d}{dz_2}
Y^{\nu}_*(v,z_2)\right)\delta
(z_1/z_2)-hz_2^{-1}Y^{\nu}_*(v,z_2)\frac
{\partial}{\partial z_1}\delta(z_1/z_2),
\end{equation}
or equivalently,
\begin{equation}\label{6.13}
[L^{\nu}(m),Y^{\nu}_*(v,z)]=\left(z^{n+1}\frac{d}{dz}+h(n+1)z^n\right)
Y^{\nu}_*(v,z)
\end{equation}
for $m\in {\Bbb Z}.$\ \ \ \ \ $\Box$
\end{prop}

Note that any $a\in \hat{L}$ is a lowest weight vector with 
weight $\frac{1}{2}\langle \bar{a}',\bar{a}'\rangle $. Thus by (\ref{7.11}),
\begin{equation}
[L^{\nu}(z_1),Y^{\nu}_*(a,z_2)]=z_2^{-1}\left(\frac{d}{dz_2}Y^{\nu}(a,z_2)\right)\delta
(z_1/z_2)-\frac{1}{2}\langle \bar{a}',\bar{a}'\rangle z_2^{-1}Y^{\nu}_*(a,z_2)\frac
{\partial}{\partial z_1}\delta(z_1/z_2).
\end{equation}
Applying (\ref{6.13}) to $v=\alpha(-1),$ which is a lowest weight
vector with weight $h=1$ for $\alpha\in {\bf h}_*^{\perp},$ and
equating the coefficients of $z^{-n-1},$ we obtain
\begin{equation}
[L^{\nu}(m),\alpha(n)]=-n\alpha(m+n)
\end{equation}
for $n\in \frac{1}{p}{\Bbb Z}.$   
 
Now in (\ref{7.8}), taking $v=\omega$ and noting that
$$L(n)\omega=0, \ \ {\rm if}\ n>0\ \ \mbox{and}\ \ n\ne 2,$$ 
\begin{equation}
L(0)\omega=2\omega,
\end{equation}
$$L(2)\omega=\frac{1}{2}\mbox{dim}\,{\bf h}^{\perp}_*,$$
we have:
\begin{equation}
\begin{array}{c}
{\displaystyle[L^{\nu}(z_1),L^{\nu}(z_2)]=z_2^{-1}\left(\frac{d}{dz_2}L^{\nu}(z_2)\right)
\delta(z_1/z_2)-2z_2^{-1}L^{\nu}(z_2)\frac{\partial}{\partial z_1}
\delta(z_1/z_2)}\\
{\displaystyle-\frac{1}{12}(\mbox{dim}\,{\bf h}^{\perp}_*)z^{-1}_2
\left(\frac{\partial}{\partial z_1}\right)^3\delta(z_1/z_2).}
\end{array}
\end{equation}
Equating the coefficients of $z^{-m-2}_1z^{-n-2}_2$, we obtain:
\begin{equation}
[L^{\nu}(m),L^{\nu}(n)]=(m-n)L^{\nu}(m+n)+\frac{1}{12}(m^3-m)
\mbox{dim}\,{\bf h}^{\perp}_*\delta_{m+n,0}
\end{equation}
for $m,n\in {\Bbb Z}$. We conclude:
\begin{prop}\label{ep6.4}
The operators $L^{\nu}(n), n\in{\Bbb Z}$ and 
$({\rm dim}\,{\bf h}^{\perp}_*)I$ span a copy of the Virasoro algebra,
and the operators $L^{\nu}(n)$ provide a representation of the Virasoro
algebra on $V^T_L$ with
\begin{equation}
\begin{array}{ll}
L_n\mapsto L^{\nu}(n)\ \, & \mbox{for}\ \ n\in{\Bbb Z}\\
\,\ c\mapsto {\rm dim}\,{\bf h}^{\perp}_*.\ &\ \ \ \ \ \ \ \ \ \ \ \ \ \  \ \Box
\end{array}
\end{equation}
\end{prop}

Finally, we justify the weight gradation of $V_L^T$  introduced in
Section 4 by using the definition of weight given by (\ref{wt}). For this
purpose we need to study $L^{\nu}(z)$ in more detail. It is easy to
see from (\ref{9.16}) that
\begin{equation}\label{6.20}
e^{\Delta_z}\alpha(-1)\beta(-1)=\alpha(-1)\beta(-1)+
\left(2c_{110}+\sum_{i=1}^{p-1}c_{11i}(\omega_p^{ki}+\omega_p^{-ki})\right)\langle\alpha,\beta\rangle
\end{equation}
for $k\in\{0,...,p-1\}$  and  
$\alpha\in({\bf h}_*^{\perp})_{(k)}, \beta\in({\bf
h}_*^{\perp})_{(-k)}.$ The numbers $c_{11r}$ for $r\in{0.,...,p-1}$
are defined in (\ref{cmn}) and equal to: 
$$c_{110}=-\frac{1}{2p^2}\sum_{i=1}^{p-1}\frac{\omega_p^{i}}
{(1-\omega_p^i)^2}$$
$$c_{11r}=\frac{1}{2p^2}\frac{\omega_p^{-r}}{(1-\omega_p^{-r})^2}\ \ \
{\rm for}\ \ \ r\ne 0.$$ 
Define constants
\begin{equation}
c_k=\sum_{r=1}^{p-1}\frac{\omega_p^{kr}}{(1-\omega_p^r)^2}.
\end{equation}
Then 
\begin{equation}\label{ck}
c_0=-(p-1)(p-5)/12,\ \ \ \ c_1=-(p^2-1)/12.
\end{equation}
Using (\ref{ck}) and the recursive formula
\begin{equation}
c_k=c_{k-1}+\frac{p+1}{2}-k+1,
\end{equation}
we find that
\begin{equation}
c_k=(k-1)(p+1-k)/2-(p^2-1)/12.
\end{equation}
Therefore,
\begin{equation}\label{6.21}
2c_{110}+\sum_{i=1}^{p-1}c_{11i}(\omega_p^{ki}+\omega_p^{-ki})=
\frac{c_{k+1}+c_{p-k+1}-2c_1}{2p^2}=\frac{k(p-k)}{2p^2}.
\end{equation}
By (\ref{6.20}) and  (\ref{6.21}) we have an
explicit expression for the vertex operator\\
$Y_*^{\nu}(\alpha(-1)\beta(-1),z)$:
\begin{equation}
Y_*^{\nu}(\alpha(-1)\beta(-1),z)=\mbox{$\circ\atop\circ$}\alpha(z)\beta(z)
\mbox{$\circ\atop\circ$}+\frac{k(p-k)}{2p^2}\langle\alpha,\beta\rangle.
\end{equation}
Recall the canonical quadratic element $\omega$ (\ref{omega}). Then by
the last formula,
\begin{equation}
L^{\nu}(z)=\frac{1}{2}\sum_{i=1}^d\mbox{$\circ\atop\circ$}\beta_i(z)\beta_i(z)\mbox{$\circ\atop\circ$}+
\frac{1}{4p^2}\sum_{k=1}^{p-1}k(p-k){\rm dim}\,({\bf h}_*^{\perp})_{(k)}.
\end{equation} 
It is obvious now that
\begin{equation}\label{z6.0}
{\displaystyle L^{\nu}(0)1=\frac{1}{4p^2}\sum_{k=1}^{p-1}k(p-k){\rm dim}\,
({\bf h}_*^{\perp})_{(k)}1}
\end{equation}
for $1\in S[\nu]$ (cf. (\ref{9.6})), and therefore the weight defined
in (\ref{gra}) is exactly the weight defined by the
$L^{\nu}(0)$-eigenvalue (\ref{wt}). 
Similarly, the weight of $t=1\otimes t\in V_L^T$ for $t\in T_{\alpha}$ 
($\alpha\in {\bf h}_{(0)}$) is $\frac{1}{2}\langle\alpha',\alpha'\rangle
+\frac{1}{4p^2}\sum_{k=1}^{p-1}k(p-k){\rm dim}\,({\bf h}_*^{\perp})_{(k)}.$
\begin{rem}\label{r6.1}\hspace*{-0.2 cm}{\bf : }{\rm The weight of $1$
in (\ref{z6.0}) is closely related to the second Bernoulli polynomial, 
defined by $B_2(x)=x^2-x+1/6.$  In fact, one can check that wt\,1
is equal to $-\frac{1}{4}\sum_{k=0}^{p-1}B_2(k/p)\,{\rm dim}\,  
({\bf h}_*^{\perp})_{(k)}+{\rm dim}\,{\bf h}_*^{\perp}/24$ (cf. [M] and 
[DM2]).}
\end{rem}

\section{Twisted modules for $V_L$}
\setcounter{equation}{0}

In this section we recall certain notions of vertex algebra and of 
twisted module.
We also present a family of $\hat\nu$-twisted modules $V_L^T$ for the
vertex algebra $V_L$ in the case in which $L$ is even. 

First we recall the definition of vertex algebra used in [DL2]. (This is 
different from the original definition in [B].)
A $vertex$  $algebra$ is a ${\Bbb Z}$-graded vector space
\begin{equation}
V=\coprod_{n\in{{\Bbb Z}}}V_n; \ \ \ \mbox{for}\ \ \ v\in V_n,\ \
n=\mbox{wt}\,v;
\end{equation}
equipped with a linear map
\begin{equation}\label{a3.30}
\begin{array}{l}
V \to (\mbox{End}\,V)[[z,z^{-1}]]\\
v\mapsto Y(v,z)=\displaystyle{\sum_{n\in{{\Bbb Z}}}v_nz^{-n-1}}\ \ \ \  (v_n\in
\mbox{End}\,V)
\end{array}
\end{equation}
and with two distinguished vectors ${\bf 1}\in V_0,$ $\omega\in V_2$
satisfying the following conditions for $u, v \in V$:
\begin{eqnarray}
& &u_nv=0\ \ \ \ \ \mbox{for}\ \  n\ \ \mbox{sufficiently large};\\
& &Y({\bf 1},z)=1;\\
& &Y(v,z){\bf 1}\in V[[z]]\ \ \ \mbox{and}\ \ \ \lim_{z\to
0}Y(v,z){\bf 1}=v;
\end{eqnarray}
\begin{equation}\label{ajac}
\begin{array}{c}
\displaystyle{z^{-1}_0\delta\left(\frac{z_1-z_2}{z_0}\right)
Y(u,z_1)Y(v,z_2)-z^{-1}_0\delta\left(\frac{z_2-z_1}{-z_0}\right)
Y(v,z_2)Y(u,z_1)}\\
\displaystyle{=z_2^{-1}\delta
\left(\frac{z_1-z_0}{z_2}\right)
Y(Y(u,z_0)v,z_2)};
\end{array}
\end{equation}
\begin{equation}
[L(m),L(n)]=(m-n)L(m+n)+\frac{1}{12}(m^3-m)\delta_{m+n,0}(\mbox{rank}\,V)
\end{equation}
for $m, n\in {{\Bbb Z}},$ where
\begin{equation}
L(n)=\omega_{n+1}\ \ \ \mbox{for}\ \ \ n\in{{\Bbb Z}}, \ \ \
\mbox{i.e.},\ \ \ Y(\omega,z)=\sum_{n\in{{\Bbb Z}}}L(n)z^{-n-2}
\end{equation}
and
\begin{eqnarray} 
& &\mbox{rank}\,V\in {{\Bbb Q}};\\
& &L(0)v=nv=(\mbox{wt}\,v)v \ \ \ \mbox{for}\ \ \ v\in V_n\
(n\in{{\Bbb Z}}); \label{a3.40}\\
& &\frac{d}{dz}Y(v,z)=Y(L(-1)v,z).\label{a3.41}
\end{eqnarray}
This completes the definition. We denote the vertex 
algebra just defined by $(V,Y,{\bf 1},\omega)$ 
(or briefly, by $V$). The series $Y(v,z)$ are called {\it vertex operators}.

An {\it automorphism} $g$ of the vertex algebra $V$ is a linear automorphism of $V$ preserving ${\bf 1}$ and $\omega$  such that the actions of $g$ 
and $Y(v,z)$ on $V$ are compatible in the sense that
\begin{equation}\label{aaut}
gY(v,z)g^{-1}=Y(gv,z)
\end{equation}
for $v\in V.$ Then $gV_n\subset V_n$ for $n\in{\Bbb Z}$ and $V$ is a direct sum of 
the eigenspaces of $g:$
\begin{equation}\label{adec}
V=\coprod_{j\in {\Bbb Z}/N{\Bbb Z}}V^j,
\end{equation}
where $N$ is the order of $g$ and $V^j=\{v\in V|gv=\omega_N^jv\}.$
 
We next recall the notion of $g$-twisted module (see [FFR] and [D]; this notion
records the properties obtained in [L1], Sec. 3.3 of [FLM2], [L2]). 
Let $(V,Y,{\bf 1},\omega)$ be a vertex 
algebra and let $g$ be an automorphism of $V$ of order $N.$ A $g$-twisted 
$module$ $M$ for $(Y,V,{\bf 1},\omega)$ is a ${\Bbb Q}$-graded
vector space
\begin{equation}\label{a6.62}
M=\coprod_{n\in{{\Bbb Q}}}M_n; \ \ \ \mbox{for}
\ \ \ w\in M_n,\ \ n=\mbox{wt}\,w;
\end{equation}
equipped with a linear map
\begin{equation}
\begin{array}{l}
V\to (\mbox{End}\,M)[[z^{1/N},z^{-1/N}]]\label{map}\\
v\mapsto\displaystyle{ Y(v,z)=\sum_{n\in{\frac{1}{N}{\Bbb Z}}}v_nz^{-n-1}\ \ \ (v_n\in
\mbox{End}\,M)}
\end{array}
\end{equation}
satisfying the following conditions for $u,v\in V$ and
$w\in M$ and $j\in{\Bbb Z}$:
\begin{eqnarray}
& &Y(v,z)=\sum_{n\in j/N+{\Bbb Z}}v_nz^{-n-1}\ \ \ \ {\rm for}\ \ v\in V^j;\label{1/2}\\
& &u_nw=0\ \ \ \mbox{for}\ \ \ n\ \ \ \mbox{sufficiently\ large};\\
& &Y({\bf 1},z)=1;
\end{eqnarray}
\begin{equation}\label{ajacm}
\begin{array}{c}
\displaystyle{z^{-1}_0\delta\left(\frac{z_1-z_2}{z_0}\right)
Y(u,z_1)Y(v,z_2)-z^{-1}_0\delta\left(\frac{z_2-z_1}{-z_0}\right)
Y(v,z_2)Y(u,z_1)}\\
\displaystyle{=z_2^{-1}\frac{1}{N}\sum_{j\in {\Bbb Z}/N{\Bbb Z}}
\delta\left(\omega_N^j\frac{(z_1-z_0)^{1/N}}{z_2^{1/N}}\right)Y(Y(g^ju,z_0)v,z_2)}
\end{array}
\end{equation}
where $Y(g^ju,z_0)$ is an operator on $V$;
\begin{equation}
[L(m),L(n)]=(m-n)L(m+n)+\frac{1}{12}(m^3-m)\delta_{m+n,0}(\mbox{rank}\,V)
\end{equation}
for $m, n\in {{\Bbb Z}},$ where
\begin{eqnarray}
& &L(n)=\omega_{n+1}\ \ \ \mbox{for}\ \ \ n\in{{\Bbb Z}}, \ \ \
\mbox{i.e.},\ \ \ Y(\omega,z)=\sum_{n\in{{\Bbb Z}}}L(n)z^{-n-2};\label{6.70}\\
& &L(0)w=nw=(\mbox{wt}\,w)w \ \ \ \mbox{for}\ \ \ w\in M_n\
(n\in{{\Bbb Q}});\label{6.71}\\
& &\frac{d}{dz}Y(v,z)=Y(L(-1)v,z).\label{6.72}
\end{eqnarray}
This completes the definition. We denote this module by $(M,Y)$ 
(or briefly, by $M$).

We now work in the setting of Remarks \ref{rem1.2}, \ref{erem4.4}
and \ref{erem4.2}. 
In particular, $L$ is a nondegenerate even lattice; 
$\nu$ is an isometry of $L$ 
such that $\nu^p=1;$ ${\bf h}=L\otimes_{\Bbb Z}{\Bbb C};$ ${\bf h}_*=0;$
$M(1)$ and $S[\nu]$ are canonical irreducible modules for $\hat{\bf h}$ 
for $\hat{\bf h}[\nu]$ respectively;  $\hat L$ and $\hat L_{\nu}$
are two central extensions of $L$ by the finite cyclic group $\langle
\kappa|\kappa^q=1\rangle, $ with commutator maps $c_0$ and $c_0^{\nu},$ 
respectively, given by (\ref{L2.2}) and (\ref{add1}); 
$c(\alpha,\beta)=(-1)^{\langle\alpha,\beta\rangle}$
and 
$c_{\nu}(\alpha,\beta)=\prod_{s=0}^{p-1}
(-\omega_p^s)^{\langle\nu^s\alpha,\beta\rangle}$ for 
$\alpha,\beta\in L;$ $\hat \nu$ is a automorphism of both $\hat L$ and 
$\hat L_{\nu}$  lifting $\nu$ such that $\hat \nu a=a$ if $\nu\bar a=\bar a$
for $a\in \hat L$ or $\hat L_{\nu}$ and such that $\hat \nu^p=1;$
$V_L=M(1)\otimes {\Bbb C}\{L\};$ ${\bf 1}=\iota(1);$ $\omega=\frac{1}{2}\sum_{r=1}^{d}\beta_r(-1)^2$
where $\{\beta_1,...,\beta_d\}$ is an orthonormal basis of ${\bf h};$
$Y(\cdot,z)=Y_*(\cdot,z)$ is the linear map in (\ref{3.41}) (see also Remark 
3.2):
\begin{equation}\label{e7.8}
\begin{array}{lcr}
V_L&\to& (\mbox{End}\,V_L)[[z,z^{-1}]]\hspace*{3.6 cm} \\
v&\mapsto& Y(v,z)=\displaystyle{\sum_{n\in{\Bbb Z}}v_nz^{-n-1}\ \ \ (v_n\in
\mbox{End}\,V_L)};
\end{array}
\end{equation}
$T$ is a $\hat L_{\nu}$-module as given in Remark \ref{erem4.4}; 
$V_L^T=S[\nu]\otimes T$ (see (\ref{9.8}));  $Y_{\nu}(\cdot,z)=Y_*^{\nu}(\cdot,z)$ is the linear map in (\ref{9.18}) (cf. Remark \ref{erem4.2}):
\begin{equation}\label{7.9}
\begin{array}{lll}
V_L&\to& (\mbox{End}\,V_L^T)[[z^{1/p},z^{-1/p}]] \\
v&\mapsto& Y_{\nu}(v,z)=\displaystyle{\sum_{n\in\frac{1}{p}{\Bbb Z}}v_n^{\nu}z^{-n-1}\ \ \ (v_n^{\nu}\in
\mbox{End}\,V_L^T)}.
\end{array}
\end{equation}

It is proved in [FLM3] (see also [B]) that 
$(V_L,Y,{\bf 1},\omega)$
is a vertex algebra of rank equal to $l={\rm rank}\, L.$  
{}From (\ref{e3.10}) and (\ref{e3.15}) we see that $\hat \nu$ fixes 
${\bf 1}$ and $\omega:$
\begin{equation}
\hat\nu {\bf 1}={\bf 1},\ \ \hat\nu \omega=\omega.
\end{equation}
Also,
\begin{equation}
\hat\nu Y(v,z)\hat\nu^{-1}=Y(\hat\nu v,z) \ \ v\in V_L
\end{equation}
 {}from (\ref{e3.48}), that is, $\hat\nu$ is an automorphism of the
vertex algebra $V_L,$ and it has period $p.$  

\begin{th}\label{x} [L1], [FLM2], [L2] Let $L$ be an even lattice as in Remark \ref{rem1.2} and
let $T$ be
an $\hat L_{\nu}$-module as in Remark \ref{erem4.4}. Then the space 
$(V_L^T,Y_{\nu})$ is a $\hat \nu$-twisted $V_L$-module. Moreover,
$V_L^T$ is irreducible if and only if $T$ is 
an irreducible $\hat L_{\nu}$-module.
\end{th}

\noindent{\bf Proof\ } Recall (\ref{9.25}) with $\psi=1,$ 
Corollary \ref{co1} and Propositions 
\ref{prop9} and \ref{ep6.4} with ${\bf h}_*=0.$ We  need only prove
(\ref{1/2}). Let $u=u^*\otimes \iota(a)$ for $u^*\in M(1)$ and
$a\in \hat L.$ Let
$$v=\sum_{r=0}^{p-1}\omega_p^{-rj}\hat\nu^{r}u\in V_L^{j}.$$ 
Then (\ref{1/2}) follows from the fact that 
$Y_{\nu}(u,z)\in ({\rm End}\,V_L^T)[[z^{1/p},z^{-1/p}]]$ (see
Remark \ref{erem4.2}) and (\ref{e4.53}).\ \ $\Box$

\end{document}